\documentclass[reqno]{amsart}
\usepackage{hyperref}
\usepackage{amssymb,amsmath,amsthm, enumerate}
\usepackage{xcolor}
\usepackage[pdftex]{graphics}
\newcommand{\R}{\mathbb{R}}
\newcommand{\Z}{\mathbb{Z}}
\newcommand{\N}{\mathbb{N}}

\def\R{{\mathbb{R}}}

\def\N{{\mathbb{N}}}
\def\Z{{\mathbb{Z}}}
\def\F{{\mathbb{F}}}

\def\grad{\nabla}

\providecommand{\keywords}[1]
{
  \small	
  \textbf{\textit{Keywords---}} #1
}

\begin{document}

\title[\hfilneg Existence and Multiplicity of Solutions Cooperative System ]{Existence and Multiplicity of Solutions for a Cooperative Elliptic System Using Morse Theory}



\author{Leandro Rec\^{o}va}
\address{California State Polytechnic University, Pomona - USA.}
\curraddr{}
\email{llrecova@cpp.edu}

\author{Adolfo Rumbos}
\address{ Pomona College, Claremont, California - USA.}
\curraddr{610 N. College Avenue, 91711 Claremont, CA}
\email{arumbos@pomona.edu}



\begin{abstract}
 In this paper, we study the existence of nontrivial solutions of the Dirichlet boundary value problem for the following elliptic system:
 \begin{equation}
  \left\{  \begin{aligned}
      -\Delta u & = au + bv + f(x,u,v);&\quad\mbox{ for }x\in\Omega,\\
       -\Delta v & = bu + cv + g(x,u,v),&\quad\mbox{ for }x\in\Omega,  \\ 
       u&=v=0,&\quad\mbox{ on }\partial\Omega, 
    \end{aligned}
    \right.
    \label{mainsys0}
\end{equation}
  for $x\in\Omega$, where $\Omega\subset\R^{N}$ is an open and connected bounded set with a  smooth boundary $\partial\Omega$, with $N\geqslant 3,$ $u,v:\overline{\Omega}\rightarrow\R$,  $a,b,c\in\R,$  and $f,g:\overline{\Omega}\times\R^2\rightarrow\R$ are continuous functions with $f(x,0,0)=0$ and $g(x,0,0) = 0$, and with super-quadratic, but sub-critical growth in the last two variables. 
 We prove that
the boundary value problem (\ref{mainsys0}) has at
least two nontrivial solutions for the case in
which the eigenvalues of the matrix 
$\displaystyle \textbf{M} 
= \begin{pmatrix}
    a & b \\ b & c
\end{pmatrix}$ are higher than the first eigenvalue
of the Laplacian over $\Omega$ with Dirichlet 
boundary conditions; $u = v= 0$ on 
$\partial\Omega$.
We use variational methods and infinite-dimensional Morse theory to obtain the multiplicity result.

\end{abstract}
\numberwithin{equation}{section}
\newtheorem{theorem}{Theorem}[section]
\newtheorem{lemma}[theorem]{Lemma}
\newtheorem{definition}[theorem]{Definition}
\newtheorem{proposition}[theorem]{Proposition}
\newtheorem{prop}[theorem]{Proposition}
\newtheorem{corollary}[theorem]{Corollary}
\newtheorem{remark}[theorem]{Remark}
\allowdisplaybreaks
\maketitle

\section{Introduction}\label{secint}
Motivated by the works of Rabinowitz, Su, and Wang \cite{RabSuWang} and Costa and 
 Magalh\~{a}es 
 in \cite{CostaMag2}, we study existence and multiplicity of the Dirichlet boundary value problem (BVP) for the following elliptic system 
\begin{equation}\label{mainsys}
\begin{cases}
  -\Delta u = au + bv + f(x,u,v), 
        & \mbox{ for } x\in\Omega;\\
  -\Delta v = bu + cv + g(x,u,v), 
        & \mbox{ for } 	x\in\Omega;\\
  u =v=0,
        &\mbox{ on }\partial\Omega,
  \end{cases}
\end{equation}
where $\Omega\subset\R^{N}$ is an open and connected bounded set with a  smooth boundary $\partial\Omega$, with $N\geqslant 3,$ $u,v:\overline{\Omega}\rightarrow\R$,  $a,b,c\in\R,$  and $f,g:\overline{\Omega}\times\R^2\rightarrow\R$ are continuous functions with $f(x,0,0)=0$ and $g(x,0,0) = 0$. 

We assume that the BVP (\ref{mainsys}) has a variational structure; that is, there exists a continuous function $F:\overline{\Omega}\times\R^2\rightarrow\R$ that is $C^1$ in the variables $u,v$ with
\begin{equation}\label{gradFDef0}
    F(x,0,0) =0,
    		\quad\mbox{ for all } x\in\Omega,
\end{equation}
and 
\begin{equation}
    \nabla F(x,s,t) = \begin{pmatrix} f(x,s,t) \\ g(x,s,t)\end{pmatrix},\quad\mbox{ for }x\in\Omega\mbox{ and }(s,t)\in\R^2,
    \label{gradFDef}
\end{equation}
where the gradient in (\ref{gradFDef}) is taken 
with respect to the variables $s$ and $t$. 

We will assume that the functions $f,g,$ and $F$ satisfy the following conditions:
\begin{itemize}
    \item [($H_1$)] $f,g\in C(\overline{\Omega}\times\R^2,\R)$,  $f(x,0,0)=0$ and $g(x,0,0) = 0$. 
    \item [($H_2$)] The functions $f$ and $g$ satisfy a super-quadratic, but sub-critical growth condition given by
   \begin{equation}
    \begin{aligned}
    |f(x,s,t)|&\leqslant C_{1}(|s|^{p-1}+|t|^{p-1});\\
      |g(x,s,t)|&\leqslant C_{2}(|s|^{p-1}+|t|^{p-1}),
        \end{aligned}
        \nonumber 
\end{equation}
for $x\in\Omega$ and all $(s,t)\in\mathbb{R}^2$,
where $p>2$, $p < 2^{*}=2N/(N-2)$, for $N\geqslant 3$, and $C_1, C_2$  positive constants. 
    
    \item [($H_3$)] $F(x,z)\geqslant 0$,  for  
 $x\in\Omega$ and $z\in\R^2.$

 \item [$(H_4)$] $\displaystyle \frac{\partial f}{\partial s}(x,s,t),\displaystyle \frac{\partial f}{\partial t}(x,s,t),\displaystyle \frac{\partial g}{\partial s}(x,s,t),$ and $\displaystyle \frac{\partial g}{\partial t}(x,s,t)$ exist and are continuous on $\overline{\Omega}\times\R^2.$
    
    \item [$(AR)_\theta$]  There exist $R>0$ and $\theta > 2$ such that 
\begin{equation}
    0 < \theta F(x,s,t) < \textbf{U}\cdot \nabla F(x,s,t),\quad\mbox{ for }|\textbf{U}| \geqslant R,
\label{arcond}
\end{equation}
where $\textbf{U}=\begin{pmatrix}s\\ t\end{pmatrix}$ and $|\textbf{U}|$ denotes the Euclidean norm of $U$. This condition is known in the literature as an Ambrosetti-Rabinowitz type condition (see \cite{AmbRab}). 
\end{itemize}

Before presenting the main results, we discuss some notation. 

Denote by $(\lambda_{k})$ the sequence of eigenvalues 
$$0< \lambda_{1} < \lambda_{2}\leqslant\lambda_{3}\leqslant\ldots\leqslant \lambda_{k}\ldots$$ 
of the $N$--dimensional Laplacian $-\Delta$ over $\Omega$,  with zero Dirichlet boundary conditions on $\partial\Omega$.  We note that $\lambda_{k}\rightarrow +\infty$ as $k\rightarrow +\infty$. 

We will denote by $\varphi_{k}$ the eigenfunction associated with the eigenvalue $\lambda_{k}$, and assume that $\|\varphi_{k}\|=1$, for $k=1,2,\ldots$, where
$\|\cdot\|$ denotes the norm in the Sobolev space 
$H^1_0(\Omega)$ given by 

\begin{equation}\label{ARnorm0}
    \|u\|=\left(\int_{\Omega}|\nabla u|^{2}\,dx\right)^{\frac{1}{2}},\quad\mbox{ for all }u\in H_{0}^{1}(\Omega). 
\end{equation}
The family of eigenfunctions $(\varphi_{k})$ forms
a complete orthonormal system for $H^1_0(\Omega)$.

Let $\textbf{M}$ be the matrix of the constant coefficients of the linear part of the system (\ref{mainsys}); namely, 
\begin{equation}\label{defamatrix}
\textbf{M}=
\begin{pmatrix}a&b\\b&c\end{pmatrix}.
\end{equation}
Denote by  $\mu_{1},\mu_{2}$ the (real) eigenvalues of $\textbf{M}$ and assume that $\mu_{1}\leqslant\mu_{2}.$  

We will use the variational approach to obtain 
solutions of the BVP in (\ref{mainsys}).  Accordingly,
we define a functional 
$J\colon H_{0}^{1}(\Omega)\times H_{0}^{1}(\Omega)\to\R$ as follows:  

\noindent
Put 
$\displaystyle 
z = \begin{pmatrix}
    u\\
    v
\end{pmatrix},\ $
where $u,v\in H_{0}^{1}(\Omega)$, and define 
$Q\colon 
H_{0}^{1}(\Omega)\times H_{0}^{1}(\Omega)\to\R$ by
\begin{equation}\label{ARq0}
    Q(z)=\frac{1}{2}\|u\|^2+\frac{1}{2}\|v\|^2
    -\frac{1}{2}\int_{\Omega}\textbf{M}z\cdot z\,dx,
\end{equation}
where $\textbf{M}$ is the matrix given in 
(\ref{defamatrix}), $\|\cdot\|$ is the Sobolev
norm in $ H_{0}^{1}(\Omega)$ given in 
(\ref{ARnorm0}), and $\textbf{M}z\cdot z$ denotes 
the Euclidean inner product of $\textbf{M}z$ and
$z$ in $\R^2$.

\noindent
Next, put 
\begin{equation}\label{ARfuncJ}
    J(z) = Q(z)-\int_{\Omega}F(x,z)\,dx,
        \quad\mbox{ for } 
   z = \begin{pmatrix}
   		u\\
   		v
   \end{pmatrix},
   \ u, v \in  H_{0}^{1}(\Omega),
\end{equation}
where $Q$ is defined in (\ref{ARq0}).  The assumptions
on $f$ and $g$  in ($H_1$) and ($H_2$) can be used
to show that the functional $J$ defined in 
(\ref{ARfuncJ}) is a $C^1$ functional with 
Fr\'{e}chet derivative  given by
\begin{equation}\label{ARDerJ05}
	\langle J^{\prime}(z),\Lambda\rangle 
	= \langle z,\Lambda\rangle_X 
    		- \int_{\Omega} 
      \textbf{M} z \cdot\Lambda\,dx 
	- \int_{\Omega}\nabla F(x,z)\cdot\Lambda\,dx,
\end{equation}
for $z = (u,v)\in X$, $\Lambda=(\varphi,\psi)\in X$,
where we have set 
$X = H_{0}^{1}(\Omega)\times H_{0}^{1}(\Omega) $  and
$$
	 \langle z,\Lambda\rangle_X 
		= \int_\Omega \nabla u\cdot \nabla \varphi\ dx 
		+ \int_\Omega \nabla v\cdot \nabla \psi\ dx.
$$

In view of (\ref{ARDerJ05}), 
$z\in H_{0}^{1}(\Omega)\times H_{0}^{1}(\Omega)$ is 
a critical point of $J$ if and only if 
$$
	 \langle z,\Lambda\rangle_X 
    		- \int_{\Omega} 
      \textbf{M} z \cdot\Lambda\,dx 
	- \int_{\Omega}\nabla F(x,z)\cdot\Lambda\,dx = 0,
    \quad\mbox{ for all } \Lambda\in 
    H_{0}^{1}(\Omega)\times H_{0}^{1}(\Omega). 
$$
Thus, critical points of the functional $J$ given in 
(\ref{ARfuncJ}) are weak solutions of the BVP in 
(\ref{mainsys}).

The main results of this paper are stated below. 

 \begin{theorem}\label{maintheo444}
 Let $\Omega$ be a bounded, open subset of $\R^{N}$, with $N\geqslant 3$.  For the Hilbert space  $X=H_{0}^{1}(\Omega)\times H_{0}^{1}(\Omega)$, consider 
 the decomposition $X=X^{-}\oplus X^{+}$, with 
 $$X^{-}=X_{m}^{-}\times X_{m}^{-}\quad\mbox{ and }\quad  X^{+}=(X^{-})^{\perp},$$ 
 where $X_m^{-}$ is the finite-dimensional space given by 
 \begin{equation}
    X_{m}^{-}=\bigoplus_{j\leqslant m}\ker(-\Delta - \lambda_{j}I).
    \nonumber 
\end{equation}
 Assume ($H_1$), ($H_2$), ($H_3$) and $(AR)_\theta$ hold true. 
 Assume also that there exists $m\geqslant 1$ such 
 that 
 $\lambda_m<\mu_1\leqslant\mu_2<\lambda_{m+1}$, where $\mu_1$ and $\mu_2$ are the eigenvalues of the matrix $\textbf{M}$ given in (\ref{defamatrix}).
Then, there exists $\delta_1>0$  such that, for
$\mu_1$ and $\mu_2$ satisfying
$$\lambda_{m}<\lambda_{m+1}-\delta_1 < \mu_1\leqslant\mu_2<\lambda_{m+1},$$ 
the functional $J$ defined in (\ref{ARfuncJ}) has at least one nontrivial critical point $z_o$. 
\end{theorem}

In the last section of this article, we prove a  second multiplicity result by computing the critical groups of the associated energy functional at the origin and at infinity, respectively.   Using an argument with a long exact sequence of a particular topological pair and condition $(H_4)$ , we obtain the following result.

  \begin{theorem}    Let $\Omega$ be a bounded and connected domain in $\R^{N}$, with $N\geqslant 3$.
  Assume that ($H_1$), ($H_2$), ($H_3$), $(H_4)$, and $(AR)_\theta$ hold true . Moreover, let $z_o$ be the critical point of $J$ given by Theorem \ref{maintheo444} . 
 \begin{itemize}
     \item [(i)] If $z_o$ is a nondegenerate critical point of $J$ with Morse index $\mu_o$, then the system (\ref{mainsys}) has at least two nontrivial solutions provided that $q_o\ne \mu_o-1,$ where $q_o=\dim X^{-}.$
     \item [(ii)] If $z_o$ is a degenerate critical point of $J$ with Morse index $\mu_o$ and nullity $\nu_0$, then the system (\ref{mainsys}) has at least two nontrivial solutions provided that 
     \begin{equation} 
     q_o+1 \not\in[\mu_o,\mu_o+\nu_o].
     \nonumber 
     \end{equation}
     where $q_o=\dim X^{-}.$ 
 \end{itemize}
\label{maintheo33}
\end{theorem}

There is an extensive literature dealing with the problem of existence and multiplicity of solutions for elliptic systems.  We discuss here some results that motivated the work in this paper. 

In \cite{CostaMag2}, Costa and  Magalh\~{a}es studied the following BVP:
\begin{equation}
    \left\{\begin{aligned}
    -\Delta u & = \lambda u \pm \delta v + f(x,u,v), &\mbox{ in }\Omega; \\
    -\Delta v & = \delta u + \gamma v + g(x,u,v), &\mbox{ in }\Omega;\\
    u&=v=0, &\mbox{ on }\partial\Omega,
    \end{aligned}
    \label{cmagsys}
    \tag{$P_\pm$}
    \right.
\end{equation}
where $\Omega$ is a bounded domain in $\R^{N}$
with smooth boundary, and $\lambda,\gamma,\delta$ are real parameters. The system $(P_+)$ is called a {\it cooperative system} and $(P_{-})$ is a {\it noncooperative system}. The authors analyzed both types of systems in the presence of a nonquadradicity condition at infinity to treat the nonresonant and resonant cases when the nonlinearity has sub-quadratic growth at infinity. They used minimax techniques to establish the existence of solutions for both non-resonant and resonant cases. 

In \cite{ZouLiLiu}, Zou, Li, and Liu studied the following BVP for the cooperative elliptic system
\begin{equation}
  \left\{  \begin{aligned}
      -\Delta u & = \lambda u + \delta v + f(x,u,v),&\quad\mbox{ in }\Omega;\\
       -\Delta v & = \delta u + \gamma v + g(x,u,v),&\quad\mbox{ in }\Omega;\\
       &u=v=0,&\quad\mbox{ on }\partial\Omega,
    \end{aligned}
    \right.
    \label{sys133}
\end{equation}
where $\Omega$ is a bounded domain in $\R^{N}(N\geqslant 3)$, and $\lambda,\delta,\gamma$ are real parameters, $f,g\in C^{1}(\overline{\Omega}\times\R^2,\R)$ and there exists some function $F\in C^{2}(\overline{\Omega}\times\R^2,\R)$ such that $\nabla F=(f,g)$. Their main interest was to study the resonance cases at the origin and at the infinity, namely,  
$$\lim_{U\rightarrow 0}\frac{\nabla F(U)}{|U|}=0\quad\mbox{ and }\quad\lim_{|U|\rightarrow\infty}\frac{\nabla F(U)}{|U|}=0,$$
respectively, 
where $|\cdot|$ denotes the usual Euclidean norm in $\R^2$ and  $U=(u,v)\in\R^2.$ The authors of \cite{ZouLiLiu} establish some conditions on the nonlinearity  that allowed them to use infinite-dimensional Morse theory to prove the existence and multiplicity of solutions of problem (\ref{sys133}).

 For a general survey on boundary value problems for nonlinear elliptic systems, we refer to the works of de Figueiredo in \cite{Dj2} and \cite{df4} where the author presents different classes of nonlinear elliptic systems and approaches for proving existence and multiplicity of solutions. 

Few of the works in the literature of elliptic systems deal with the use of infinite-dimensional Morse theory to study  existence and multiplicity of solutions of such problems. The main contribution of this paper is to use infinite-dimensional Morse theory to study an elliptic system whose scalar version is a problem considered by Rabinowitz, Su, and Wang in \cite{RabSuWang}.

This paper is organized as follows:  In Section \ref{prelim} we present some notation and results that will be used throughout this paper.  In Section \ref{auxFunc}, we discuss the auxiliary functional $Q$ defined in (\ref{ARq0}). We show that $Q$ satisfies the Palais-Smale condition and compute its critical groups at the origin. 
In Section \ref{pssection}, we proceed to prove that the energy functional $J$ associated with problem (\ref{mainsys}) satisfies the Palais-Smale condition. The first existence result is proved in Section \ref{linksection} where we show that the functional $J$ satisfies the hypotheses of the linking theorem of Rabinowitz \cite{Rab78}.   Finally, using properties of a long-exact sequence of homology groups of a particular topological pair, we prove the existence of a second nontrivial solution for problem (\ref{mainsys}) in Section \ref{maintheosec}.

\section{Preliminaries}\label{prelim}
Let $\mathcal{H}$ denote a Hilbert space with inner product $\langle \cdot,\cdot\rangle_{\mathcal{H}}$. The Cartesian product $X=\mathcal{H}\times\mathcal{H}$ with inner product $\langle\cdot,\cdot\rangle_{X}:\mathcal{H}\times \mathcal{H}\rightarrow\R$ defined by
\begin{equation}
\langle (u_1,v_1),(u_2,v_2)\rangle_{X} = \langle u_1,v_1\rangle + \langle u_2,v_2\rangle,
\label{hilb1}
\end{equation}
for all $(u_1,v_1),(u_2,v_2)\in X$, is a Hilbert space. 
 In this paper, we are interested in the cases in which $\mathcal{H}=H_{0}^{1}(\Omega)$, or $\mathcal{H}=L^{2}(\Omega)$. We describe the former case next.

Let $\Omega\subset\R^{N}$ be an open and bounded subset of $\R^{N}$, for $N\geqslant 3$, and let $H_{0}^{1}(\Omega)$ denote the Sobolev space 
obtained through completion of $C^\infty_c(\Omega)$ with respect to the norm 
\begin{equation}
    \|u\|=\left(\int_{\Omega}|\nabla u|^{2}\,dx\right)^{\frac{1}{2}},\quad\mbox{ for all }u\in H_{0}^{1}(\Omega), 
    \label{norm1}
\end{equation}
where $C^\infty_c(\Omega)$ denotes the space of 
real-valued $C^\infty$ functions in $\R^N$ that have
compact support contained in $\Omega$.

The norm in (\ref{norm1}) is generated by the inner product 
\begin{equation}
    \langle u,\varphi\rangle = \int_{\Omega}\nabla u\cdot\nabla\varphi\,dx,\quad\mbox{ for all }u,\varphi\in X.
    \label{inner1}
\end{equation}
Denote $H_{0}^{1}(\Omega)\times H_{0}^{1}(\Omega)$ 
by $X$, with the induced inner product $\langle\cdot,\cdot\rangle_{X}:X\times X\rightarrow\R$ defined by 
\begin{equation}
    \langle z,\Phi\rangle_{X} = \langle u,\varphi\rangle + \langle v,\psi\rangle,
    \label{inner2}
\end{equation}
for $z=(u,v)\in X$ and $\Phi=(\varphi,\psi)\in X,$ and norm $\|\cdot\|_{X}:X\rightarrow\R$ defined by 
\begin{equation}
    \|z\|_{X}^2=\|u\|^2 + \|v\|^2,
    \quad\mbox{ for } z=(u,v)\in X. 
    \label{norm2}
\end{equation}

Similarly, the Cartesian product $Y=L^{2}(\Omega)\times L^2(\Omega)$ 
is also a Hilbert space with inner product defined by
\begin{equation}
    \langle z,\Phi\rangle_{Y} = \langle u,\varphi\rangle_{L^2} + \langle v,\psi\rangle_{L^2},
    \label{inner3}
\end{equation}
for all $z=(u,v),\Phi=(\varphi,\psi)\in Y$, 
by virtue of (\ref{hilb1}), 
where $\displaystyle
\langle u,v\rangle_{L^2}=\int_{\Omega}uv\,dx,$
for $u,v\in L^2(\Omega)$,
is the inner product in $L^2(\Omega)$.  We also
obtain the norm $\| \cdot\|_{Y}$ given by
\begin{equation}\label{YnormDefn05}
	\|z\|_Y^2 = \|u\|_{L^2(\Omega)}^2 + \|v\|_{L^2(\Omega)}^2 ,
        \quad\mbox{ for all } z=(u,v)\in Y.
\end{equation}

Next, consider the BVP for the linear elliptic system
\begin{equation}
  \left\{  \begin{aligned}
      -\Delta u & = au + bv, &\quad\mbox{ in }\Omega;\\
       -\Delta v & = bu + cv  ,&\quad\mbox{ in }\Omega,\\
       u=0&,v=0, &\quad\mbox{ on }\partial\Omega.
         \end{aligned}
    \right.
    \label{sys0}
\end{equation}

We will represent the system (\ref{sys0}) in matrix  notation 
\begin{equation}
    -\overrightarrow{\Delta} z=\textbf{M}z,
    \label{sys1}
\end{equation}
where 
\begin{equation}
    -\overrightarrow{\Delta} z=\begin{pmatrix}-\Delta u\\  -\Delta v\end{pmatrix}\quad\mbox{ and }\quad\textbf{M}=\begin{pmatrix} a&b\\b&c\end{pmatrix},
    \label{eq1}
\end{equation}
for $z=(u,v)\in X$, with $a,b,c\in\R$.

Let $\sigma(\textbf{M})=\{\mu_1,\mu_2\}$ denote
the spectrum of the matrix $\textbf{M}$ and $\textbf{v}_1,\textbf{v}_2\in \R^2$ be  eigenvectors corresponding to $\mu_1$ and $\mu_2$, respectively, such that 
\begin{equation}\label{OrthoNormalBasisM}
    \langle \textbf{v}_1, \textbf{v}_2\rangle_{\R^2} = 0
    \quad\mbox{ and }\quad |\textbf{v}_1|=|\textbf{v}_2|=1
\end{equation}
where $|\cdot|$ denotes the Euclidean norm in $\R^2$. 

The eigenvalues $\mu_1$ and $\mu_2$ of $\textbf{M}$
can be obtained as the minimum and the maximum values, respectively,  of the quadratic form $\langle \textbf{M}z,z\rangle_{\R^2}$ restricted to the unit sphere in $\R^2$. We also have that
\begin{equation}\label{pest342}
    \mu_1 |z|^2 \leqslant   \textbf{M}z\cdot z \leqslant \mu_2|z|^2,\quad\mbox{ for all }z\in \R^2.
\end{equation}
 
 We will associate an energy functional $Q:X\rightarrow\R$ with the linear problem (\ref{sys1}) given by
\begin{equation}
    Q(z)=\frac{1}{2}\|u\|^2+\frac{1}{2}\|v\|^2-\frac{a}{2}\int_{\Omega}u^2\,dx-\frac{b}{2}\int_{\Omega}u^2\,dx-\frac{b}{2}\int_{\Omega}v^2\,dx-\frac{c}{2}\int_{\Omega}v^2\,dx,
    \label{q1}
\end{equation}
for all $z=(u,v)\in X.$ 

By virtue of the definition of the norm given by (\ref{norm2}), we can write (\ref{q1}) in a more compact form as 
\begin{equation} \label{q2}
    Q(z) = \frac{1}{2}\|z\|_{X}^2-\frac{1}{2}\int_{\Omega}\textbf{M}z\cdot z\,dx,
        \quad\mbox{ for } z\in X.
\end{equation}

According to \cite[Proposition $1.2$]{CostaMag2}, there exist an isomorphism $i:X\rightarrow X$, mutually orthogonal subspaces $V,N,W$ of $X$, with $\dim N < \infty$, and a constant $\mu > 0$ such that the space $X$ can be decomposed as $X=V\oplus N\oplus W$; 
so that, setting $q=Q\circ i,$ we have
\begin{equation}
\begin{aligned}
    q(z) & \leqslant -\mu\|z\|_{X}^2, &\mbox{ for all }z\in V,\\
     q(z) &   =0, &\mbox{ for all }z\in N,\\
      q(z) & \geqslant \mu\|z\|_{X}^2, &\mbox{ for all }z\in W.
      \end{aligned}
\label{q3}
\end{equation}
Furthermore, $N\not=\{0\}$ if an only if the
matrix $\textbf{M}-\lambda_j I$, where $I$ denotes the 
$2\times 2$ identity matrix, is singular for some $\lambda_j$. 

The isomorphism $i:X\rightarrow X$ used in  \cite[Proposition $1.2$]{CostaMag2}  is such that 
\begin{equation}
    i(\varphi,\psi) = \varphi \textbf{v}_1 + \psi \textbf{v}_2,
        \quad \hbox{ for } (\varphi,\psi)\in X,
    \label{iso1}
\end{equation}
where $\textbf{v}_1,\textbf{v}_2$ are eigenvectors of the matrix $\textbf{M}$ associated with the eigenvalues $\mu_1,\mu_2,$ respectively, and 
satisfying (\ref{OrthoNormalBasisM}).

To see how the subspaces $V$, $N$ and $W$ come 
about, for $\xi\in\R$, define a quadratic form
$r_\xi\colon H^1_0 (\Omega) \to\R$ by 
\begin{equation}\label{QuadFormr1}
     r_\xi(u) = \|u\|^2-\xi\|u\|_{L^2}^2,
        \quad\hbox{ for } u\in H^1_0 (\Omega).
\end{equation}
Using the completeness of the orthonormal 
system $(\varphi_j)$, we can show that the 
quadratic form given in (\ref{QuadFormr1}) is
negative definite on the space
\begin{equation}\label{HminusDfn05}
     H_{\xi}^{-} 
        = \bigoplus_{\lambda_j < \xi}
        \ker(-\Delta -\lambda_{j} I),
\end{equation}
and positive definite on the space 
\begin{equation}\label{HplusDfn05}
     H_{\xi}^{+} 
        = \bigoplus_{\lambda_j > \xi}
        \ker(-\Delta -\lambda_{j} I).
\end{equation}
If $\xi$ is an eigenvalue, $\lambda$, of 
$(-\Delta, H^1_0 (\Omega))$, set 
\begin{equation}\label{HzeroDfn05}
     H_{\xi}^{0} 
        = \bigoplus_{\lambda_j = \lambda}
        \ker(-\Delta -\lambda_{j} I);
\end{equation}
otherwise, set
\begin{equation}\label{HzeroDfn10}
     H_{\xi}^{0} = \{ 0\}.
\end{equation}
We therefore obtain the decomposition
$ H_{0}^{1}(\Omega) 
=  H_{\xi}^{-}\oplus H_{\xi}^{0}\oplus H_{\xi}^{+}$,
where the subspaces $ H_{\xi}^{-}$ and  
$H_{\xi}^{+}$ are given in (\ref{HminusDfn05}) 
and (\ref{HplusDfn05}), respectively, and 
$H_{\xi}^{0}$ is given by (\ref{HzeroDfn05}) in 
the case in which $\xi$ is an eigenvalue, $\lambda$ of $(-\Delta, H^1_0 (\Omega))$, or by
(\ref{HzeroDfn10}) otherwise. 

Using the notation introduced above for the 
quadratic functional $r_\xi$ in (\ref{QuadFormr1}),  with $\mu_1$ and $\mu_2$ in
place of $\xi$, we define the subspaces $V$, $N$ 
and $W$ as follows:
\begin{equation}\label{decx3}
V=H_{\mu_1}^{-}\times H_{\mu_2}^{-}, \quad
 N=H_{\mu_1}^{0}\times H_{\mu_2}^{0}, 
 \quad\mbox{ and }\quad 
 W = H_{\mu_1}^{+}\times H_{\mu_2}^{+}.
\end{equation}
The assertion in (\ref{q3}) now follows from the
fact that 
\begin{equation}\label{qz005}
    q(\varphi, \psi) = \frac{1}{2}r_{\mu_1}(\varphi) + \frac{1}{2}r_{\mu_2}(\psi), \quad\mbox{ for }(\varphi, \psi)\in X,
\end{equation}
where $r_{\mu_1}$ and $r_{\mu_2}$ are given in 
(\ref{QuadFormr1}), with $\xi$ replaced by 
$\mu_1$ and $\mu_2$, respectively, and $V$, $N$
and $W$ are given in (\ref{decx3}).

We note that the subspaces $V$ and $N$ defined in
(\ref{decx3}) are finite dimensional vector 
spaces.

In subsequent sections in this paper we will need the following decomposition of the Hilbert space $X$:
\begin{equation}
X=X^{-}\oplus X^{+},
\label{decx1}
\end{equation}
where 
\begin{equation}\label{decx2}
X^{-} = V\oplus N\quad\mbox{ and } \quad
X^{+} = W,
\end{equation}
with $V$, $N$ and $W$ as given in (\ref{decx3}).
Since, $V$ and $N$ are finite dimensional, it 
follows from the definition of $X^{-}$ in 
(\ref{decx2}) that $\dim X^{-}<\infty$.

We will represent the system (\ref{mainsys}) in a compact form as
\begin{equation}
    -\overrightarrow{\Delta }z = \textbf{M}z + \nabla F(x,z);\quad\mbox{ in }\Omega,\quad z=0\mbox{ on }\partial\Omega,
    \nonumber
\end{equation}
where $F:\overline{\Omega}\times \R^2\rightarrow\R$ is a continuous function 
that is differentiable in the variable $z=(u,v)$ 
with
$\nabla F(x,z)=(f(x,u,v),g(x,u,v))$,
where the gradient is taken with respect
to the variables $u$ and $v$.  Furthermore, 
the functions 
$f,g:\overline{\Omega}\times\R^2\rightarrow\R$ given in (\ref{mainsys}) are continuous and satisfy the super-quadratic, but subcritical growth conditions given in ($H_2$); namely,
\begin{equation}
    \begin{aligned}
    |f(x,s,t)|&\leqslant C_{1}(|s|^{p-1}+|t|^{p-1});\\
      |g(x,s,t)|&\leqslant C_{2}(|s|^{p-1}+|t|^{p-1}),
        \end{aligned}
        \label{subcg}
\end{equation}
for $x\in\Omega$ and $(s,t)\in\mathbb{R}^2$, where 
$p > 2$ and $p<2^{*}=2N/(N-2)$, $N\geqslant 3$, and $C_{1}$ and $C_{2}$ are positive constants. 

Let $J:X\rightarrow\R$ be the energy functional associated with problem (\ref{mainsys}) given by 
\begin{equation}\label{func0}
    J(z) = Q(z)-\int_{\Omega}F(x,z)\,dx,
        \quad\mbox{ for } z\in X,
\end{equation}
where $Q$ is defined in (\ref{q2}). The functional (\ref{func0}) is well-defined; this can be shown 
to be a consequence of the growth condition in  (\ref{subcg}) and the Sobolev embedding theorem. 
It can also be shown that $J\in C^1(X,\R)$ with
 Fr\'{e}chet derivative  given by  
\begin{equation}
    \langle J^{\prime}(z),\Lambda\rangle = \langle Q^{\prime}(z),\Lambda\rangle_X - \int_{\Omega}\nabla F(x,z)\cdot\Lambda\,dx,
    \quad\mbox{ for }\Lambda=(\varphi,\psi)\in X,
\label{jprimef}
\end{equation}
and $z = (u,v)\in X$, where 
\begin{equation}\label{qfrechdef}
    \langle Q^{\prime}(z),\Lambda\rangle 
    	= \langle z,\Lambda\rangle_X 
    		- \int_{\Omega} 
      \textbf{M} z\cdot\Lambda\,dx,\quad\mbox{ for }\Lambda=(\varphi,\psi)\in X,
\end{equation}
and $z = (u,v)\in X$. Consequently,
\begin{equation}\label{DerJ05}
	\langle J^{\prime}(z),\Lambda\rangle 
	= \langle z,\Lambda\rangle_X 
    		- \int_{\Omega} 
      \textbf{M} z \cdot\Lambda\,dx 
	- \int_{\Omega}\nabla F(x,z)\cdot\Lambda\,dx,
\end{equation}
for $\Lambda=(\varphi,\psi)\in X$ 
and $z = (u,v)\in X$.

As was shown in \cite{AmbRab}, or in \cite[Remark $2.13$]{Rabinowitz1}, the $(AR)_{\theta}$ condition provides a super-quadratic estimate for the function $F$. We present the version for the elliptic system (\ref{mainsys}) here for the reader's convenience. 

\begin{lemma} Let $F:\overline{\Omega}\times\R^2\rightarrow\R$ be the $C^1$ function satisfying the condition $(AR)_{\theta}$. Then,
\begin{equation}
    F(x,\textbf{U})\geqslant C_{1}|\textbf{U}|^{\theta},\quad\mbox{ for all }|\textbf{U}|\geqslant R,
    \label{arr1}
\end{equation}
for $\textbf{U}=\begin{pmatrix}u\\v\end{pmatrix}$, where $C_1$ is a positive constant.  
\label{lemmaest1}
\end{lemma} 
\begin{proof}[Proof:]
    Assume that the condition $(AR)_{\theta}$ holds 
    true.  We then have that
    \begin{equation}\label{ineq1}
        \frac{\partial F}{\partial u}u + \frac{\partial F}{\partial v}v \geqslant \theta F,\quad\mbox{ for }|\textbf{U}|\geqslant R.
    \end{equation}
The characteristic curves associated with the left-hand
side of the differential inequality in (\ref{ineq1}) 
are given by solutions of the system of ordinary differential equations 
\begin{equation*}
    \frac{du}{ds}=u
    \quad\mbox{ and }\quad
    \frac{dv}{ds}=v,\quad\mbox{ for }s\in\R,
\end{equation*}
which yields the curves parametrized by 
\begin{equation}\label{ineq3}
    u(s) = c_1e^s\quad\mbox{ and }\quad v(s)=c_2e^s,\quad\mbox{ for }s\in\R,
\end{equation}
for some constants $c_1$ and $c_2$.

Along the characteristic curves in (\ref{ineq3}), 
we get that 
\begin{equation}\label{ARineq13}
    |\textbf{U}| = \sqrt{c_1^2+c_2^2}\ e^s,
        \quad\mbox{ for } s\in\R;
\end{equation}
so that, the condition $|\textbf{U}|\geqslant R$ in 
(\ref{ineq1}) translates to 
\begin{equation}\label{ARineq6}
    e^s \geqslant \frac{R}{\sqrt{c_1^2 +c_2^2}},
        \quad\mbox{ provided that }  
            (c_1,c_2)\not=0.
\end{equation}
Taking the initial points for the characteristic curves
to lie on the unit circle
\begin{equation}\label{ARineq14}
    c_1^2 + c_2^2 =1,
\end{equation}
we obtain from (\ref{ARineq6}) that the condition 
$|\textbf{U}|\geqslant R$ is equivalent to 
$$
    e^s \geqslant R,
        \quad\mbox{ or } \quad
            s \geqslant \ln R.
$$
Hence, the differential inequality in (\ref{ineq1}) 
translate into 
\begin{equation}\label{ARineq7}
    \frac{d}{ds}F(x,u,v) \geqslant \theta F(x,u,v),
        \quad\mbox{ for } s\geqslant \ln R,
\end{equation}
along the characteristic curves.

Next, rewrite the inequality in (\ref{ARineq7}) as
$$
    \frac{d}{ds}F(x,u,v) - \theta F(x,u,v) \geqslant 0,
        \quad\mbox{ for } s\geqslant \ln R,
$$
and multiply by the integrating factor $e^{-\theta s}$
to obtain 
\begin{equation}\label{Far1}
    \frac{d}{ds}\left( 
        e^{-\theta s} F(x,u,v) 
    \right)
    \geqslant 0,
        \quad\mbox{ for } s\geqslant \ln R,
\end{equation}
along characteristic curves. 

Next, integrate the inequality in (\ref{Far1}) along
characteristic curves from $\ln R$ to s, for 
$s\geqslant\ln R$ to get that 
\begin{equation}\label{AReqn9}
    e^{-\theta s} F(x,u,v) 
    \geqslant \frac{1}{R^\theta} F(x,c_1 R , c_2 R),
            \quad\mbox{ for } s\geqslant\ln R,
\end{equation}
where we have used the expressions for $u$ and $v$ in
(\ref{ineq3}) along characteristic curves.

Next, set 
\begin{equation}\label{AReqn10}
    m_1 = \min_{\genfrac{}{}{0pt}{2}{c_1^2+c_2^2=1}{x\in\overline{\Omega}}}
        F(x,c_1 R , c_2 R).
\end{equation}
It follows from the left-most inequality in 
(\ref{arcond}) that the constant $m_1$ defined in 
(\ref{AReqn10}) is strictly positive. 

Combining (\ref{AReqn10}) and (\ref{AReqn9}) then
yields 
\begin{equation*}
    e^{-\theta s} F(x,u,v) 
    \geqslant \frac{m_1}{R^\theta} ,
            \quad\mbox{ for } s\geqslant\ln R,
\end{equation*}
along characteristic curves, which can be rewritten as
\begin{equation}\label{AReqn11}
     F(x,u,v) 
    \geqslant \frac{m_1}{R^\theta} e^{\theta s},
            \quad\mbox{ for } s\geqslant\ln R,
\end{equation}
along characteristic curves.

Finally, combine (\ref{ARineq13}) and (\ref{ARineq14})
to obtain that 
$$
    |\textbf{U}|  = e^s,
        \quad\mbox{ for } s\in\R,
$$
along characteristic curves, from which we get that
\begin{equation}\label{ARineq12}
     e^{\theta s} = |\textbf{U}|^\theta ,
        \quad\mbox{ for } s\in\R,
\end{equation}
along characteristic curves. 

Next, substitute the expression for $e^{\theta s}$ 
in (\ref{ARineq12}) into the inequality in  
(\ref{AReqn11}) to obtain 
\begin{equation}\label{AReqn15}
     F(x,u,v) 
    \geqslant 
    \frac{m_1}{R^\theta} |\textbf{U}|^\theta ,
        \quad\mbox{ for } |\textbf{U}| \geqslant  R.
\end{equation}

Setting $C_1 = \dfrac{m_1}{R^\theta}\ $ in 
(\ref{AReqn15}) yields the inequality in (\ref{arr1})
and this completes the proof of the lemma.
\end{proof}

Based on the above lemma, there exist constants $c_1,c_2>0$ such that 
\begin{equation}
    F(x,\textbf{U})\geqslant c_{1}|\textbf{U}|^{\theta} - c_2,\quad\mbox{ for all }x\in\overline{\Omega}\mbox{ and }\textbf{U}\in\R^2.
    \label{arr2}
\end{equation}

To apply results from critical point theory, and to compute the critical groups at isolated 
critical points of the associated energy functional $J$ defined in (\ref{func0}), we need to verify a compactness condition presented below. 
 
\begin{definition}{\rm We will say that $(u_{m},v_{m})$ is a PS sequence for $J$ if 
\begin{equation}
    |J(u_{m},v_{m})|\leqslant C\quad\mbox{ for all }m\quad\mbox{ and }\quad \nabla J^{\prime}(u_{m},v_{m})\rightarrow 0\quad\mbox{ as }m\rightarrow\infty,
    \label{pscond1}
\end{equation}
where $C$ is a constant. We say that a functional $J$ defined on the space 
$X = H_{0}^{1}(\Omega)\times H_{0}^{1}(\Omega)$ satisfies the Palais-Smale condition, or PS condition, if any PS-sequence $(u_{m},v_{m})\subset X$ possesses a convergent subsequence. 
\label{psdef0}}
\end{definition}

Finally,  we present some notation that will be used to define the critical groups of $J$ at isolated critical points.

Let $A,B$ be two topological spaces with $B\subset A$. Denote by $H_{q}(A,B)$ the $q$-singular relative homology group of the pair $(A,B)$ with coefficients in a field $\mathbb{F}$.  Let $c=J(u_{0},v_{0})$, where $(u_{0},v_{0})$ is an isolated critical point of $J$, and set 
$$J^{c}=\{(u,v)\in X : J(u,v)\leqslant c\}.$$ 
The $q$-th critical groups of $J$ at $(u_{0},v_0)$, with coefficients in $\mathbb{F}$, are given by
\begin{equation}
C_{q}(J,(u_{0},v_{0})) = H_{q}(J^{c}\cap U,(J^{c}\cap U)\backslash\{(u_{0},v_{0})\}),\quad q\in\Z,
\label{cgroupdef}
\end{equation}
(see \cite[Definition $4.1$, p. 32]{KC}), where $U$ is an open neighborhood of $(u_{0},v_{0})$ such that $(u_{0},v_{0})$ is the only critical point of $J$ in $U$. The critical groups of isolated critical points are well--defined and they do not depend on the choice of the neighborhood $U$. This follows from the excision property of homology theory.

\section{The Palais-Smale condition and critical groups for Q}\label{auxFunc}

Throughout this section we assume that there exists 
$m\geqslant 1$ such that 
\begin{equation}\label{NonResEqn0005}
	\lambda_m<\mu_1\leqslant\mu_2<\lambda_{m+1},
\end{equation}
where $\mu_1$ and $\mu_2$ are the eigenvalues of the matrix $\textbf{M}$ given in (\ref{defamatrix}).

The assumption in 
(\ref{NonResEqn0005}) will be used to derive some properties of the
functional $Q\colon X\to\R$ defined in (\ref{q2}) that will help us compute the critical groups of $J$ at the origin in the next  section.

Consider the following decomposition of the Hilbert space $X$ given by    
 \begin{equation} \label{decz1}
 X=X^{-}\oplus X^{+}, 
 \end{equation}
 where 
\begin{equation}\label{DfnDecompXEqn05}
    X^{-}=X_{m}^{-}\times X_{m}^{-}
        \quad\mbox{ and }\quad
            X^{+}=( X^{-})^{\perp},
\end{equation}
with 
\begin{equation}\label{dec1}
    X_{m}^{-}=\bigoplus_{j\leqslant m}\ker(-\Delta - \lambda_{j}I).
\end{equation}

We will use the above decomposition of the space $X$,
and the assumption in (\ref{NonResEqn0005}), 
to show that the functional $Q$ defined in (\ref{q2})
satisfies the Palais-Smale condition.

\begin{lemma} 
Assume that there exists $m\geqslant 1$ such that 
(\ref{NonResEqn0005}) holds true.  Then, the 
functional $Q:X\rightarrow\R$ defined in 
(\ref{q2}); namely,
$$Q(z)=\frac{1}{2}\|z\|_{X}^2-\frac{1}{2}\int_{\Omega}\textbf{M}z\cdot z\,dx,
    \quad\mbox{ for } z\in X,
$$
satisfies the Palais-Smale condition.    
\label{psQ0lemma}
\end{lemma}
 
\begin{proof}[Proof:] Let $(z_n)\subset X$ be a Palais-Smale sequence for $Q$; so that, 
\begin{equation*}
    |Q(z_n)|\leqslant C ,
        \quad\mbox{ for all } n\in\mathbb{N},
\end{equation*}
and some constant $C>0$, and  
\begin{equation}\label{pscc10}
    \left|\langle Q^{\prime}(z_n),\Phi\rangle\right|\leqslant\varepsilon_n\|\Phi\|_{X},
        \quad\mbox{ for all } \Phi\in X 
            \mbox{ and } n\in\mathbb{N},
\end{equation}
where $\varepsilon_n > 0$, for all $n$, with $\varepsilon_n\rightarrow 0$ as $n \rightarrow\infty$. 

We show that the sequence $(z_n)$ converges to $0$.  
This will establish the Palais-Smale condition for $Q$.

Write $z_n = z_{n1}+z_{n2},$ with $z_{n1}\in X^{-}$ and $z_{n2}\in X^{+}$. 

We first show that 
\begin{equation}\label{3324Eqn0015}
    \lim_{n\to\infty} \|z_{n2}\|_X = 0.
\end{equation}

To establish (\ref{3324Eqn0015}), let $\Phi = z_{n2}$
in (\ref{pscc10}) and use the definition of the 
 Fr\'{e}chet derivative of $Q$ in (\ref{qfrechdef})
 to obtain that 
\begin{equation}\label{3324Eqn0020}
    \|z_{n2}\|_X^2 
    	-\int_\Omega \textbf{M} z_{n2}\cdot z_{n2}\ dx
    		\leqslant \varepsilon_n \|z_{n2}\|_X
    			\quad\mbox{ for all } n\in\N,
\end{equation}
where $\varepsilon_n \searrow 0$ as 
$n \rightarrow\infty$, and where  we have also used the fact that 
$$
\int_{\Omega}\textbf{M}z_{n1}\cdot z_{n2}\,dx = 0,
$$
since $X^{+}$ is the orthogonal complement of $X^{-}$ in $X$.

Next, use this estimate in (\ref{pest342}) to
get 
\begin{equation}\label{3324Eqn0025}
    \int_\Omega \textbf{M} z\cdot z \ dx 
    		\leqslant\mu_2 \|z\|_Y^2 ,
    			\quad\mbox{ for all } z\in X,
\end{equation}
where $\|\cdot\|_Y$ is the norm in 
$Y= L^2(\Omega) \times L^2(\Omega)$ defined in (\ref{YnormDefn05}).

On the other hand, using the fact that the family of 
eigenfunctions, $(\varphi_{k})$, of 
$(-\Delta, H^1_0(\Omega))$ forms a complete orthogonal 
system for $H^1_0(\Omega)$, we can obtain the 
inequality 
\begin{equation}\label{PoincareIneq0015}
	\|z_2\|_X^2 \geqslant \lambda_{m+1} \|z_2\|_Y^2,
		\quad\mbox{ for all } z_2 \in X^+,
\end{equation}
Consequently, combining (\ref{3324Eqn0025}) and 
(\ref{PoincareIneq0015}),
\begin{equation}\label{PoincareIneq0020}
	\int_\Omega \textbf{M} z_{n2}\cdot z_{n2}\ dx
	\leqslant
	\frac{\mu_2}{\lambda_{m+1}}\|z_2\|_X^2,  
		\quad\mbox{ for all } z_2 \in X^+.
\end{equation}
It then follows from (\ref{3324Eqn0020}) and 
(\ref{PoincareIneq0020}) that 
\begin{equation}\label{PoincareIneq0025}
	\left( 
	1 - \frac{\mu_2}{\lambda_{m+1}}
	\right)
	\|z_2\|_X^2
	\leqslant 
	\varepsilon_n \|z_{n2}\|_X
    			\quad\mbox{ for all } n\in\N.
\end{equation}
Observe that 
$$
	1-\frac{\mu_2}{\lambda_{m+1}} > 0,
$$
in view of the assumption in (\ref{NonResEqn0005}).
Thus, it follows from (\ref{PoincareIneq0025}),
and the fact that $\varepsilon_n \searrow 0$ as 
$n \rightarrow\infty$, that 
$$
    \lim_{n\to\infty} \|z_{n2}\|_X = 0,
$$
which is the assertion in (\ref{3324Eqn0015}).

Next, we show that 
\begin{equation}\label{3324Eqn0035}
    \lim_{n\to\infty} \|z_{n1}\|_X = 0.
\end{equation}
To prove the assertion in (\ref{3324Eqn0035}), use the assumption in (\ref{pscc10}) to obtain that
\begin{equation}\label{3324Eqn0005}
    -\langle Q^{\prime}(z_n),\Phi\rangle
    \leqslant\varepsilon_n\|\Phi\|_{X},
        \quad\mbox{ for all } \Phi\in X 
            \mbox{ and } n\in\mathbb{N},
\end{equation}
where $\varepsilon_n \searrow 0$ as 
$n \rightarrow\infty$.
Thus, using the definition of the Fr\'{e}chet 
derivative of $Q$ in (\ref{qfrechdef}), and putting 
$\Phi = - z_{n1}$ in (\ref{3324Eqn0005}), 
we obtain that
\begin{equation}\label{3324Eqn0010}
    \int_{\Omega}\textbf{M}z_{n1}\cdot z_{n1}\,dx 
    	- \|z_{n1}\|_{X}^2 \leqslant 
    		\varepsilon_n \|z_{n1}\|_{X},
\end{equation}
where  we used the fact that 
$$
\int_{\Omega}\textbf{M}z_{n2}\cdot z_{n1}\,dx = 0,
$$
since $X^{+}$ is the orthogonal complement of $X^{-}$ in $X$.

Next, use the estimate in (\ref{pest342}) to
get 
\begin{equation}\label{3324Eqn0040}
    \int_\Omega \textbf{M} z\cdot z \ dx 
    		\geqslant \mu_1 \|z\|_Y^2 ,
    			\quad\mbox{ for all } z\in X.
\end{equation}
On the other hand, using the fact that the family of 
eigenfunctions, $(\varphi_{k})$, of 
$(-\Delta, H^1_0(\Omega))$ forms a complete orthogonal 
system for $H^1_0(\Omega)$, we can obtain the 
inequality 
\begin{equation}\label{PoincareIneq0045}
	\|z_1\|_X^2 \leqslant \lambda_{m} \|z_1\|_Y^2,
		\quad\mbox{ for all } z_1 \in X^-.
\end{equation}
Thus, combining (\ref{3324Eqn0040}) and 
(\ref{PoincareIneq0045}),
\begin{equation}\label{PoincareIneq0050}
	\int_\Omega \textbf{M} z_{1}\cdot z_{1}\ dx
	\geqslant
	\frac{\mu_1}{\lambda_{m}}\|z_{1}\|_X^2,  
		\quad\mbox{ for all } z_1 \in X^-.
\end{equation}
It then follows from (\ref{3324Eqn0010}) and 
(\ref{PoincareIneq0050}) that 
\begin{equation}\label{PoincareIneq0055}
	\left( \frac{\mu_1}{\lambda_{m}} -1
	\right)
	\|z_1\|_X^2
	\leqslant 
	\varepsilon_n \|z_{n1}\|_X
    			\quad\mbox{ for all } n\in\N.
\end{equation}
Observe that 
$$
	\frac{\mu_1}{\lambda_{m}} - 1  > 0,
$$
in view of the assumption in (\ref{NonResEqn0005}).
Thus, it follows from (\ref{PoincareIneq0055}),
and the fact that $\varepsilon_n \searrow 0$ as 
$n \rightarrow\infty$, that 
$$
    \lim_{n\to\infty} \|z_{n1}\|_X = 0,
$$
which is the assertion in (\ref{3324Eqn0035}).

Combining (\ref{3324Eqn0035}) and (\ref{3324Eqn0015})
we obtain that 
$$
    \lim_{n\to\infty} \|z_{n}\|_X = 0,
$$
which shows that any the Palais-Smale sequence $(z_n)$
of $Q$ converges to $0$.  This concludes the proof of the lemma. 
\end{proof}

In the next result, we show that the origin is the unique critical point of the functional $Q$ defined in (\ref{q2}).

\begin{lemma}\label{lemmaQ}
Assume that $\lambda_m<\mu_1\leqslant\mu_2<\lambda_{m+1}$ for some $m\geqslant 1$, where $\mu_1\leqslant\mu_2$ are the eigenvalues of the matrix $\textbf{M}$. Then, the origin  is an isolated critical point of the functional $Q$ defined in (\ref{q2}).
\end{lemma}

\begin{proof}[Proof:]
We will show that the origin is the only critical point of $Q$ in $X$. 

First, use the definition of the Fr\'{e}chet derivative of $Q$ in (\ref{qfrechdef}) to write
\begin{equation}\label{qd1}
\langle Q^\prime(z),\Phi\rangle = \langle z,\Phi\rangle_X - \int_{\Omega}\textbf{M}z\cdot\Phi\,dx,
    \quad\mbox{ for } z,\Phi\in X.
\end{equation}
Next, write $z=z_1+z_2$, where $z_1\in X^{-}$ and $z_2\in X^{+}$, where $X^{-}$ and $X^{+}$ are defined in (\ref{DfnDecompXEqn05}) and (\ref{dec1}).  Similarly,
write $\Phi = \xi_1 + \xi_2$, where $\xi_1\in X^-$
and $\xi_2\in X^+$. We then obtain from (\ref{qd1})
that 
\begin{equation}\label{3524Eqn0005}
\langle Q^\prime(z),\Phi\rangle 
= \langle z_1,\xi_1\rangle_X + 
\langle z_2,\xi_2\rangle_X
- \int_{\Omega}\textbf{M}z_1\cdot\xi_1\,dx
- \int_{\Omega}\textbf{M}z_2\cdot\xi_2\,dx,
\end{equation}
where  we have also used the fact that 
$$
\int_{\Omega}\textbf{M}z_1\cdot \xi_2\,dx = 0
	\quad\mbox{ and }\quad
		\int_{\Omega}\textbf{M}z_2\cdot \xi_1\,dx = 0,
$$
since $X^{+}$ is the orthogonal complement of $X^{-}$ in $X$. 

For the special choice $\Phi = -z_1 + z_2$, we obtain
from (\ref{3524Eqn0005}) that 
\begin{equation}\label{3524Eqn0010}
\langle Q^\prime(z),z_2 - z_1\rangle 
= \|z_2\|_X^2 - \int_{\Omega}\textbf{M}z_2\cdot z_2\,dx
	+ \int_{\Omega}\textbf{M}z_1\cdot z_1\,dx 
		- \|z_1\|_X^2 ,
\end{equation}
for $z=z_1+z_2 \in X^-\oplus X^+$.

Now, by virtue of (\ref{pest342}), we have
\begin{equation}\label{qd12}
\mu_1\|{z}\|_{Y}^2\leqslant\int_{\Omega}\textbf{M}z\cdot z\,dx\leqslant\mu_2\|{z}\|_{Y}^2,
    \quad\mbox{ for all } z\in X,
\end{equation}
where $\|\cdot\|_{Y}$ is the norm in $L^2(\Omega)\times L^{2}(\Omega)$ defined in (\ref{inner3}).

Next, use the fact that the family of 
eigenfunctions, $(\varphi_{k})$, of 
$(-\Delta, H^1_0(\Omega))$ forms a complete orthogonal 
system for $H^1_0(\Omega)$ to obtain the 
inequality 
\begin{equation}\label{3924Eqn0005}
	\|z_2\|_X^2 \geqslant \lambda_{m+1} \|z_2\|_Y^2,
		\quad\mbox{ for all } z_2 \in X^+,
\end{equation}
Consequently, combining (\ref{3924Eqn0005}) with the right-most
inequality in (\ref{qd12}), we obtain the inequality
\begin{equation}\label{3924Eqno0010}
	\int_\Omega \textbf{M} z_{2}\cdot z_{2}\ dx
	\leqslant
	\frac{\mu_2}{\lambda_{m+1}}\|z_2\|_X^2,  
		\quad\mbox{ for all } z_2 \in X^+.
\end{equation}
We then obtain from (\ref{3924Eqno0010}) that 
\begin{equation}\label{qd14}
    \|{z}_2\|_{X}^2-\int_{\Omega}\textbf{M}{z}_2\cdot {z}_2\,dx \geqslant \left(1-\frac{\mu_2}{\lambda_{m+1}}\right)\|\overline{z}_2\|_{X}^2, \quad\mbox{ for all }z_2\in X^{+}.
\end{equation}
where
\begin{equation}\label{RightEst05}
	1-\frac{\mu_2}{\lambda_{m+1}} > 0,
\end{equation}
in view of the assumption in (\ref{NonResEqn0005}).

Similarly, using the completeness of the orthogonal 
system of eigenfunctions  $(\varphi_{k})$ in 
$H^1_0(\Omega)$, we obtain 
\begin{equation}\label{31024Eqn0010}
	\|{z}_1\|_{X}^2\leqslant\lambda_{m}\|{z}_1\|_{Y}^2,
		\quad\mbox{ for } z_1\in X^-.
\end{equation}
Combining (\ref{31024Eqn0010}) with the left-most 
inequality in (\ref{qd12}) we then obtain that 
\begin{equation}\label{qd16}
\int_{\Omega}\textbf{M}{z}_1\cdot{z}_1\,dx \geqslant \frac{\mu_1}{\lambda_{m}}\|z_{1}\|_{X}^2,\quad\mbox{ for }z_1\in X^{-}.
\end{equation}
We then obtain from (\ref{qd16}) that 
\begin{equation}\label{31024Eqn0015}
\int_{\Omega}\textbf{M}{z}_1\cdot{z}_1 \,dx 
	- \|z_{1}\|_{X}^2\geqslant 
\left( \frac{\mu_1}{\lambda_{m}} - 1\right)\|z_{1}\|_{X}^2,\quad\mbox{ for }z_1\in X^{-},
\end{equation}
where 
\begin{equation}\label{LefttEst05}
	\frac{\mu_1}{\lambda_{m}} - 1  > 0,
\end{equation}
in view of the assumption in (\ref{NonResEqn0005}).

Substituting the estimates in (\ref{31024Eqn0015}) 
and (\ref{qd14}) into the right-hand side
of (\ref{3524Eqn0010}), we obtain from (\ref{3524Eqn0010}) that
\begin{equation}\label{31024Eqn0020}
\langle Q^\prime(z),z_2 - z_1\rangle 
\geqslant \nu \|z\|_X^2,
	\quad\mbox{ for } z=z_1+z_2 \in X^-\oplus X^+,
\end{equation}
where
\begin{equation}\label{31024Eqn0025}
	\nu = \min\left\{
				\frac{\mu_1}{\lambda_{m}} - 1,
				1-\frac{\mu_2}{\lambda_{m+1}}  
				\right\}.
\end{equation}
It follows from the definition of $\nu$ in 
(\ref{31024Eqn0025}) and the observations in 
(\ref{RightEst05}) and (\ref{LefttEst05}) that 
$\nu > 0$.

Suppose that $\overline{z}$ is a critical point of $Q$; so that,
\begin{equation}\label{qd2}
     \langle Q^\prime(\overline{z}),\Phi\rangle =0,\quad\mbox{ for all }\Phi\in X.
\end{equation}
In particular, setting 
$\Phi=-\overline{z}_1 + \overline{z}_2$ in 
(\ref{qd2}), where 
$\overline{z} = \overline{z}_1 + \overline{z}_2$, we get 
\begin{equation}\label{qd3}
     \langle 
     Q^\prime(\overline{z}),
     \overline{z}_2 - \overline{z}_1\rangle = 0. 
\end{equation}
It then follows from (\ref{qd3}) and 
(\ref{31024Eqn0020}) that 
$$
    \nu\|\overline{z}\|_X^2 \leqslant 0,
$$
from which we get that 
$\|\overline{z}\|_{X}=0$; so that,  
any critical point of $Q$ must be the origin. 
\end{proof}

 Finally, we present the computation of the critical groups of $Q$ at the origin.

 \begin{lemma}
    Assume that $\lambda_m<\mu_1\leqslant\mu_2 < \lambda_{m+1}$, for some $m\geqslant 1$.
     Let $X=X^{-}\oplus X^{+}$ be the decomposition of the Hilbert space $X$ as in (\ref{dec1}).  Then, the origin is a nondegenerate critical point of $Q$. Moreover, the critical groups of $Q$ at the origin are given by 
     \begin{equation}
         C_{k}(Q,0)\cong \delta_{k,\dim X^{-}}\F,\quad\mbox{ for }k\in\Z.
         \label{cgroupsQ} 
     \end{equation}
 \end{lemma}

\begin{proof}[Proof:] 
The Fr\'{e}chet derivative of the functional $Q$ defined in (\ref{q2}) is given by 
\begin{equation}\label{qs2}
    \langle Q^\prime(z),\Phi\rangle = \langle z_1,\varphi\rangle _{X} + \langle z_2,\psi\rangle_{X}- \int_{\Omega}\textbf{M}z\cdot \Phi\,dx,
\end{equation}
for $z=z_1 + z_2\in X^-\oplus X^+$ and 
$\Phi= \varphi + \psi  \in X^- \oplus X^+$, with
$z_1 ,\varphi\in X^-$ and $z_2, \psi\in X^+$.

It follows from the result of Lemma \ref{lemmaQ} that
the origin in $X$ is an isolated critical point of 
$Q$.  The second Fr\'{e}chet derivative of $Q$ at the origin is given by
\begin{equation}
    Q^{\prime\prime}(0)(\Phi,\Psi) = \langle \varphi,\xi_1\rangle_{X} + \langle \psi,\xi_2\rangle_{X} - \int_{\Omega}\textbf{M}\Phi\cdot \Psi\,dx,
    \label{secff1}
\end{equation}
for all $\Phi=\varphi + \psi,\Psi=\xi_1 + \xi_2\in X,$ with $\varphi,\xi_1\in X^{-}$ and $\psi,\xi_2\in X^{+}.$

Set $\Phi=\Psi=\xi_1\in X^{-}$ in (\ref{secff1}). Then, we can write

\begin{equation}
    Q^{\prime\prime}(0)(\Phi,\Phi) = \|\xi_1\|_{X}^2  - \int_{\Omega}\textbf{M}\xi_1\cdot \xi_1\,dx, 
    \label{secff2}
\end{equation}
for $\xi_1\in X^{-}$. 

By virtue of the estimate (\ref{pest342}), we obtain from (\ref{secff2}) that

\begin{equation}
    Q^{\prime\prime}(0)(\Phi,\Phi) \leqslant \|\xi_1\|_{X}^2  - \mu_1\|\xi_1\|_{Y}^2\,dx, 
    \label{secff3}
\end{equation}
for $\xi_1\in X^{-}$. 

Recall that $\|\xi_1\|_{X}^2\leqslant \lambda_{m}\|\xi_1\|_{Y}^2$, for all $\xi_1\in X^{-}$. Then, it follows from (\ref{secff3}) that

\begin{equation}
    Q^{\prime\prime}(0)(\Phi,\Phi) \leqslant \left(1-\frac{\mu_1}{\lambda_m}\right)\|\xi_1\|_{Y}^2,
    \label{secff4}
\end{equation}
for $\xi_1\in X^{-}$. 
Notice that, 
\begin{equation}
    1-\frac{\mu_1}{\lambda_m} < 0,
    \nonumber 
\end{equation}
since $\mu_1 > \lambda_{m}$. Thus, we have  
\begin{equation}
    Q^{\prime\prime}(0)(\Phi,\Phi) < 0,\quad\mbox{ for }\Phi\in X^{-}\backslash\{0\},
    \nonumber 
   \end{equation}
Similar calculations show that 
$$
    Q^{\prime\prime}(0)(\Psi,\Psi) > 0,\quad\mbox{ for }\Psi\in X^{+}\backslash\{0\}.
$$
Thus, the origin in $X$ is a non-degenerate critical point of $Q$ with Morse index $q_o = \dim X^{-}.$ Hence, it follows from  \cite[Corollary $8.3$]{MW1} that 
\begin{equation}
    C_{k}(Q,0) \cong \delta_{k,\dim X^{-}}\F,\quad\mbox{ for all }k\in\Z.
    \nonumber 
\end{equation}

\end{proof}

\section{The Palais-Smale condition for the functional J}\label{pssection} 

Throughout this section, we assume that hypotheses 
($H_1$), ($H_2$), and $(AR)_\theta$ hold true. 
We restate these conditions here for the reader's 
convenience.

\begin{itemize}
    \item [($H_1$)] $f,g\in C(\overline{\Omega}\times\R^2,\R)$,  $f(x,0,0)=0$ and $g(x,0,0) = 0$. 
    \item [($H_2$)] The functions $f$ and $g$ satisfy a subcritical growth condition given
    by
   \begin{equation}
    \begin{aligned}
    |f(x,s,t)|&\leqslant C_{1}(|s|^{p-1}+|t|^{p-1});\\
      |g(x,s,t)|&\leqslant C_{2}(|s|^{p-1}+|t|^{p-1}),
        \end{aligned}
        \nonumber 
\end{equation}
for all $(s,t)\in\R^2$,
where $p>2$, $p < 2^{*}=2N/(N-2)$, for $N\geqslant 3$, and $C_1, C_2$ positive constants.

    \item [$(AR)_\theta$]  There exist $R>0$ and $\theta > 2$ such that 
\begin{equation*}
    0 < \theta F(x,u,v) < \textbf{U}\cdot \nabla F(x,u,v),\quad\mbox{ for }|\textbf{U}| \geqslant R,
\end{equation*}
where $\textbf{U}=\begin{pmatrix}u\\v\end{pmatrix}$ and $|\textbf{U}|$ denotes the Euclidean norm of $U$. 
\end{itemize}

Under these assumptions, we show that the functional
$J$ given in (\ref{func0}) satisfies the Palais-Smale condition. 

\begin{lemma}
Assume that  ($H_1$), ($H_2$) and $(AR)_\theta$ are satisfied. Then, the functional $J$ defined in (\ref{func0}) satisfies the Palais-Smale condition. 
\label{pslemma}
\end{lemma}

\begin{proof}[Proof:] We follow an approach similar to that in Rabinowitz \cite[Theorem $5.16$]{Rabinowitz1}. Let $(z_n)\subset X$ be a Palais-Smale sequence. Then, according to the definition of Palais-Smale sequence in (\ref{pscond1}), 
\begin{equation}\label{pscon2}
|J(z_n)|\leqslant C,\quad\mbox{ for all }n\in\N,
\end{equation}
where $C>0$ is a constant, and 
\begin{equation}\label{pscond3}
    \left|\langle J^{\prime}(z_n),\Phi\rangle\right|\leqslant\varepsilon_{n}\|\Phi\|_{X},
    \quad\mbox{ for all } n\in\N
        \mbox{ and all } \Phi \in X,
\end{equation}
where
\begin{equation}\label{epsProp05}
    \varepsilon_n\searrow 0
        \ \mbox{ as }\ n\rightarrow +\infty.
\end{equation}

Set $\Phi=z_n$ into (\ref{pscond3}) to get
\begin{equation}\label{PSproofEqn0005}
    \varepsilon_{n}\| z_n\|_{X} 
    \geqslant 
    	-\langle  J^{\prime}(z_n),z_n\rangle,	
    		\quad\mbox{ for all } n\in\N.
\end{equation}
Now, by the property of the sequence $(\varepsilon_n)$ in (\ref{epsProp05}), there exists $N_1\in\N$ such that 
\begin{equation}\label{condde}
    0<\varepsilon_n < 1,\quad\mbox{ for }n\geqslant N_{1}.
\end{equation}
Combining (\ref{condde}) and 
(\ref{PSproofEqn0005}) we then obtain
\begin{equation}\label{PSproofEqn0010}
    \| z_n\|_{X} \geqslant
    	-\langle  J^{\prime}(z_n),z_n\rangle,	
    		\quad\mbox{ for } n\geqslant N_1.
\end{equation}
Choose  $\gamma\in(\frac{1}{\theta},\frac{1}{2})$ where $\theta$ is given in $(AR)_\theta$ condition. Next, multiply on both 
sides of (\ref{PSproofEqn0010}) by 
$\gamma$  to get 
\begin{equation}\label{PSproofEqn0015}
    \gamma \| z_n\|_{X} \geqslant
    	-\gamma\langle 
    	 J^{\prime}(z_n),z_n\rangle,	
    		\quad\mbox{ for } n\geqslant N_1.
\end{equation}

Using the estimates in (\ref{pscon2}) and 
(\ref{PSproofEqn0015}) we obtain 
\begin{equation}\label{PSproofEqn0020}
    C + \gamma\|z_n\|_{X}
    \geqslant 
     J(z_n) 
    	- \gamma
    	\langle  J^{\prime}(z_n),z_n\rangle,
    		\quad\mbox{ for } n\geqslant N_1.
\end{equation}

It follows from the definition of the functional $J$ in (\ref{func0}), the definition of its Fr\'{e}chet derivative
in (\ref{DerJ05}), and the inequality in (\ref{PSproofEqn0020}) that 
     \begin{equation} 
    \begin{aligned}
    C + \gamma\|z_n\|_{X}
    &\geqslant \left(\frac{1}{2}-\gamma\right)\|z_{n}\|_{X}^2-\left(\frac{1}{2}-\gamma\right)\int_{\Omega}\textbf{M}z_{n}\cdot z_{n}\,dx\\& \qquad +\int_{\Omega}\left[\gamma\nabla F(x,z_n)\cdot z_n - F(x,z_n) \right]\,dx ,
            \end{aligned}
          \label{psest3}
         \end{equation}
for all $n\geqslant N_1.$ 

For $R$ given in condition  $(AR)_{\theta}$, define
\begin{equation}\label{PScondPf0005}
	\Omega_{R,n} = 
	\{ x\in\Omega\mid |z_n(x)|\leqslant R\},
		\quad\mbox{ for } n\in\N. 
\end{equation}
It then follows from the assumption in 
$(AR)_{\theta}$ that 
\begin{equation}\label{PScondPf0010}
	0 < \theta F(x,z_n(x)) < \grad F(x,z_n(x))\cdot z_n(x),
		\quad\mbox{ for } x\in\Omega_{R,n}^c, 
\end{equation}
and $n\in\N$, 
where $\Omega_{R,n}^c$ denotes the complement of $\Omega_{R,n}$
defined in (\ref{PScondPf0005}).

Next, put 
\begin{equation}\label{PScondPf0015}
C_3 = 
\max_{\genfrac{}{}{0pt}{2}{x\in\overline{\Omega}}{w\in\R^2, |w|\leqslant R}}
\left|\gamma \grad F(x,w)\cdot w - F(x,w)\right|.
\end{equation}
We then obtain the estimate 
\begin{equation}\label{PScondPf0020}
\begin{array}{rcl}
\displaystyle
\int_{\Omega}\left[\gamma\nabla F(\cdot,z_n)\cdot z_n - F(\cdot,z_n) \right]\,dx
& \geqslant & \displaystyle
\int_{\Omega_{R,n}^c}\left[\gamma\nabla F(\cdot,z_n)\cdot z_n - F(\cdot,z_n) \right]\,dx \\
& & \qquad\qquad - C_3|\Omega|,
\end{array}
\end{equation}
for all $n\in\N$, 
where $|\Omega|$ denotes the Lebesgue measure of $\Omega$.

Then, since $\gamma > \dfrac{1}{\theta}$, we obtain from 
the estimates in (\ref{PScondPf0010}) and (\ref{PScondPf0020})
that 
\begin{equation}\label{PScondPf0025}
\int_{\Omega}\left[\gamma\nabla F(\cdot,z_n)\cdot z_n - F(\cdot,z_n) \right]\,dx \geqslant
(\gamma\theta -1) 
\int_{\Omega_{R,n}^c} F(x,z_n)\ dx - C_3|\Omega|,
\end{equation}
for all $n\in\N$, where $C_3$ is given in 
(\ref{PScondPf0015}).

Combining the estimates in (\ref{psest3}) and 
(\ref{PScondPf0025}), we obtain 
\begin{equation}
     \begin{aligned}
    C + \gamma\|z_n\|_{X} 
    &\geqslant \left(\frac{1}{2}-\gamma\right)\|z_{n}\|_{X}^2-\left(\frac{1}{2}-\gamma\right)\int_{\Omega}\textbf{M}z_{n}\cdot z_{n}\,dx \\ & \qquad + \left(\gamma\theta -1\right)\int_{\Omega_{R,n}^c} F(x,z_n)\,dx - C_3|\Omega|,
    \end{aligned}
    \label{psest4}
\end{equation}
for  $n\geqslant N_1$, 
where $C_3$ is the non-negative constant given in 
(\ref{PScondPf0015}).

Next,  write 
\begin{equation}
    \int_{\Omega_{R,n}^c} F(x,z_n)\,dx =  \int_{\Omega} F(x,z_n)\,dx -  \int_{\Omega_{R,n}} F(x,z_n)\,dx,
    \label{interm1}
\end{equation}
for all $n$, and estimate the second integral on
the right-hand side of (\ref{interm1}) as follows
\begin{equation}\label{PScondPf0030}
	\left| 
		\int_{\Omega_{R,n}} F(x,z_n)\,dx
	\right|
	\leqslant 
	C_4 |\Omega|,
		\quad\mbox{ for all } n\in\N,
\end{equation}
where 
$$
C_4 = 
\max_{\genfrac{}{}{0pt}{2}{x\in\overline{\Omega}}{w\in\R^2, |w|\leqslant R}} \left| F(x,w)\right|.
$$
It then follows from (\ref{interm1}) that 
\begin{equation}\label{PScondPf0035}
\int_{\Omega_{R,n}^c} F(x,z_n)\,dx \geqslant
\int_{\Omega} F(x,z_n)\,dx - C_4|\Omega|.
\end{equation}

Combining the estimates in (\ref{psest4}) and 
(\ref{PScondPf0035}) then yields
\begin{equation}\label{ARpsest5}
     \begin{aligned}
    C + \gamma\|z_n\|_{X} 
    &\geqslant \left(\frac{1}{2}-\gamma\right)\|z_{n}\|_{X}^2-\left(\frac{1}{2}-\gamma\right)\int_{\Omega}\textbf{M}z_{n}\cdot z_{n}\,dx \\ & \qquad + \left(\gamma\theta -1\right)\int_{\Omega} F(x,z_n)\,dx - C_6,
    \end{aligned}
\end{equation}
for $n\geqslant N_1$, and some positive 
constant $C_5$.

Next, use the estimate in (\ref{arr2}) to 
obtain from (\ref{ARpsest5}) the estimate
\begin{equation}\label{ARpsest6}
     \begin{aligned}
    C + \gamma\|z_n\|_{X} 
    &\geqslant \left(\frac{1}{2}-\gamma\right)\|z_{n}\|_{X}^2-\left(\frac{1}{2}-\gamma\right)\int_{\Omega}\textbf{M}z_{n}\cdot z_{n}\,dx \\ & \qquad + c_1\left(\gamma\theta -1\right)\int_{\Omega} |z_n|^\theta\,dx - C_6,
    \end{aligned}
\end{equation}
for $n\geqslant N_1$, where $C_6$ is a positive 
constant.

Using the right-most estimate in  (\ref{pest342}), the estimate in 
(\ref{ARpsest6}) can be rewritten as 
\begin{equation}
     \begin{aligned}
    C + \gamma\|z_n\|_{X} 
    &\geqslant \left(\frac{1}{2}-\gamma\right)\|z_{n}\|_{X}^2- \left(\frac{1}{2}-\gamma\right)\mu_2\|z_n\|_{Y}^2 \\ & \qquad + c_1\left(\gamma\theta -1\right)\int_{\Omega}|z_n|^{\theta}\,dx - C_6,
    \end{aligned}
    \label{psest5}
\end{equation}
for $n\geqslant N_1$.

Next, we find estimates for $\| z_n\|_{Y}$ in terms of the $L^\theta$ norm of $z_n$, for 
$n\in\mathbb{N}$.  To do so, write 
$z_n = (u_n,v_n)$, for $n\in\mathbb{N}$, where 
$ u_n , v_n \in  H_{0}^{1}(\Omega)$ for 
$n\in\mathbb{N}$.

Set $q_1=\theta/2$ and $q_2=q_1/(q_1-1)$, its conjugate exponent. Then, by virtue of H\"{o}lder's inequality with exponents $q_1$ and $q_2$, we have
\begin{equation}
    \|u_n\|_{L^2(\Omega)}\leqslant |\Omega|^{1/2q_2}\|u_n\|_{L^\theta(\Omega)},
        \quad\mbox{ for } n\in \mathbb{N}.
    \label{psest6}
\end{equation}
Next, apply Young's inequality (see \cite[Section $B.2$, page $622$]{EV}),
\begin{equation*}
    ab\leqslant \varepsilon a^{q_1}+C(\varepsilon)b^{q_2},
\end{equation*}
for $a,b>0$, $\varepsilon>0$, where $C(\varepsilon)$ is a positive constant and $C(\varepsilon)\rightarrow\infty$ as $\varepsilon\rightarrow 0$, to the right-hand
side of (\ref{psest6}) with $a = \|u_n\|_{L^\theta(\Omega)}$ and
$b = |\Omega|^{1/2q_2}$, and $q_1$ and $q_2$
as given above, to obtain
\begin{equation}
    \|u_n\|_{L^2(\Omega)}\leqslant \varepsilon \|u_n\|_{L^{\theta}(\Omega)}^{\theta/2} + C(\varepsilon)|\Omega|^{1/2},
        \quad\mbox{ for } n\in\mathbb{N}.
    \label{psest7}
\end{equation}
Similar calculations involving $v_n$ yield
\begin{equation}
     \|v_n\|_{L^2(\Omega)}\leqslant \varepsilon \|v_n\|_{L^{\theta}(\Omega)}^{\theta/2} + C(\varepsilon)|\Omega|^{1/2},
        \quad\mbox{ for } n\in\mathbb{N}.
    \label{psest8}
\end{equation}
Using the fact that $(a+b)^2\leqslant 2(a^2+b^2),$ we obtain from (\ref{psest7}) and (\ref{psest8}) that
$$
    \begin{aligned}
        \|z_n\|_{Y}^2 & = \|u_n\|_{L^2(\Omega)}^2 + \|v_n\|_{L^2(\Omega)}^2 \\ 
        & \leqslant  2\varepsilon\left[\|u_n\|_{L^\theta(\Omega)}^\theta+\|v_n\|_{L^\theta(\Omega)}^\theta\right] + 4C(\varepsilon)^2|\Omega|,
    \end{aligned}
$$
for $n\in\mathbb{N}$; so that,
\begin{equation}\label{psest9}
	\|z_n\|_{Y}^2 \leqslant  
	4\varepsilon\ \|z_n\|_{L^\theta(\Omega)}^\theta  
		+ 4C(\varepsilon)^2|\Omega|,
    			\quad\mbox{ for } n\in\mathbb{N}.
\end{equation}

Substitute the estimate in (\ref{psest9}) into (\ref{psest5}) to obtain the estimate
\begin{equation}
     \begin{aligned}
    C + \gamma\|z_n\|_{X} 
    &\geqslant \left(\frac{1}{2}-\gamma\right)\|z_{n}\|_{X}^2   + \left[c_1\left(\gamma\theta -1\right)-\left(2-4\gamma\right)\mu_2\varepsilon\right]
    \|z_n\|_{L^\theta(\Omega)}^\theta  \\
    &\qquad\qquad  -2\left(1-2\gamma\right)(C(\varepsilon))^2|\Omega| -  C_5,
    \end{aligned}
    \label{psest15}
\end{equation}
for $n\geqslant N_1$.

Next, choose $\varepsilon>0$ small enough so
that
$$
    \varepsilon < \frac{c_1(\gamma\theta -1)}{2(1-2\gamma)\mu_2}
$$
to obtain from (\ref{psest15}) that 
\begin{equation}\label{ARpsest16}
    C + \gamma\|z_n\|_{X} 
    \geqslant \left(\frac{1}{2}-\gamma\right)\|z_{n}\|_{X}^2   
     -2\left(1-2\gamma\right)(C(\varepsilon))^2|\Omega| -  C_5,  
        \mbox{ for } n\geqslant N_1.
\end{equation}
It then follows from (\ref{ARpsest16}) that the sequence $(z_n)$ is bounded in $X$.

Hence, since the Palais-Smale sequence $(z_n)$ is bounded in $X$, we can follow the argument found
in \cite[Proposition B$.35$]{Rabinowitz1} to show 
that the sequence $(z_n)$ has a convergent 
subsequence.  The argument relies on the Sobolev
embedding theorem and the subcritical growth behavior of $f$ and $g$ in (\ref{subcg}), since we 
are assuming that $2 < p < 2^{*}$.
\end{proof}

\section{First Multiplicity Result}\label{linksection}
Throughout this section, in addition to assuming ($H_1$), ($H_2$), and $(AR)_\theta$, we
will also assume ($H_3$):
\begin{itemize}
    \item [($H_3$)] $F(x,z)\geqslant 0$,  for  $x\in\Omega$ and $z\in\R^2.$
\end{itemize}

We will prove the existence of a nontrivial solution of the system (\ref{mainsys}) by using the Linking Theorem of Rabinowitz \cite{Rab78} as presented in Willem \cite[Theorem $1.12$]{Willem1}.

\begin{theorem}[Linking Theorem]
Let $X=X^{-}\oplus X^{+}$ be a Banach space with $\dim X^{-} < \infty.$ Let $\rho > r > 0$ and let $z_2\in X^{+}$ be such that $\|z_2\|=r$. Define
\begin{itemize} 

    \item [$M$] $ : =\{z=z_1+t z_2 : \|z\|\leqslant\rho, t\geqslant 0, z_1\in X^{-}\}$,
    \item [$M_0$] $ : = \{z=z_1+t z_2: z_1\in X^{-},\|z\|=\rho\mbox{ and } t\geqslant 0 \mbox{ or } \|z\|\leqslant \rho\mbox{ and }t=0\}$, 
    \item [$N$] $:=\{z\in X^{+} : \|z\|=r\}.$
 \end{itemize}
Let $J\in C^{1}(X,\R)$ be such that 
\begin{equation}
    b_o : =\inf_{z\in N}J(z) > a_o:=\max_{z\in M_0}J(z).
    \label{bcond1}
\end{equation}
If $J$ satisfies the Palais-Smale condition at level $c$ with 
\begin{equation} 
	 c : = \inf_{\gamma\in\Gamma}\max_{z\in M}J(\gamma(z)),
  \nonumber 
\end{equation}
where
$$
	\Gamma : = \{\gamma\in C(M,X): \gamma|_{M_0}=id\},
$$
then $c$ is a critical value of $J.$
\label{linktheo}
\end{theorem}

In what follows, we will show that the functional $J$ satisfies the hypotheses of Theorem \ref{linktheo}.

\begin{lemma}
Assume that  ($H_1$), ($H_2$), ($H_3$) and $(AR)_\theta$ are satisfied, and that the 
functional $J\colon X\to\R$ is as defined in 
 (\ref{func0}), where $Q$ is given in (\ref{ARq0}). 
Let $X=X^{-}\oplus X^{+}$ be the decomposition of the Hilbert 
 space $X$ given in (\ref{decz1}), (\ref{DfnDecompXEqn05}) and (\ref{dec1}). Moreover, assume that $\lambda_m<\mu_1\leqslant\mu_2<\lambda_{m+1}$, where $\mu_1$ and
 $\mu_2$ are the eigenvalues of the matrix 
 $\textbf{M}$ given in (\ref{eq1}).  Then, there exist  $a_1, r_1>0$ such that 
\begin{equation}
 \inf_{z\in N} J(z) \geqslant a_1,
\label{cond1link}
\end{equation}
where $N=\{z\in X^{+}:\|z\|_{X}=r_1\}.$
\label{lemmalink1}
\end{lemma}

\begin{proof}[Proof:] 
 Let $z=(u,v)\in X^{+}$; so that, $u,v\in (X_{m}^{-})^{\perp}$, where $X_{m}^{-}$ is given in
 (\ref{dec1}). Then, it follows from Poincar\'{e}'s inequality that
\begin{equation}
    \lambda_{m+1}\|u\|_{L^2}^2 \leqslant \|u\|^2\quad\mbox{ and }  \quad \lambda_{m+1}\|v\|_{L^2}^2 \leqslant \|v\|^2.
    \label{clink1}
\end{equation}

Next, use the definition of the functional 
$J$ given in  (\ref{func0}) and that of the quadratic
functional $Q$ given in (\ref{q2}) to obtain 
\begin{equation}\label{LinkingEqn0005}
	J(z) = \frac{1}{2}\|z\|_{X}^2-\frac{1}{2}\int_{\Omega}\textbf{M}z\cdot z\,dx -\int_{\Omega} F (x,z)\,dx,
        \quad\mbox{ for } z\in X.
\end{equation}

It follows from the hypotheses in  $(H_1)$ and $(H_2)$, and the definition of the gradient of 
$F$ in (\ref{gradFDef}), that 
\begin{equation}
    |F(x,s.t)|\leqslant C_1\left(|s|^{p}+|t|^{p}\right),\quad\mbox{ for }x\in\Omega, (s,t)\in \R^2,
    \label{BigFEst}
\end{equation}
where $C_1$ is a positive constant. 

Using the Sobolev embedding theorem, we obtain 
from (\ref{BigFEst})  that 
\begin{equation}
    \left|\int_{\Omega}F(x,z)\,dx\right|\leqslant C\|z\|_{X}^{p},\quad\mbox{ for }x\in\Omega, z\in X,
    \label{BigFest1}
\end{equation}
where $C$ is a positive constant. 

Substitute (\ref{BigFest1}) and (\ref{pest342}) into (\ref{LinkingEqn0005}) to get 
\begin{equation}\label{LinkingEqn0006}
	J(z) \geqslant \frac{1}{2}\|z\|_{X}^2-\frac{\mu_2}{2}\left(\|u\|_{L^2}^2+\|v\|_{L^2}^2\right) -C\|z\|_{X}^{p},
        \quad\mbox{ for } z\in X.
\end{equation}
 
Then, using the inequalities in (\ref{clink1}),
we get from (\ref{LinkingEqn0006}) that 
\begin{equation}\label{clink2}
	J(z)\geqslant 
		 \frac{1}{2}\left(1-\frac{\mu_2}{\lambda_{m+1}}  - 2C\|z\|_{X}^{p-2}\right)\|z\|_{X}^2,
\end{equation}
for all $z=(u,v)$, with 
$u,v\in (X_{m}^{-})^{\perp}$.  Note that this is
the juncture in the argument in which we use the
assumption that $p>2$ in ($H_2$).

Put 
\begin{equation}\label{deflittler1}
   r = \left[\frac{1}{4C}\left(1-\frac{\mu_2}{\lambda_{m+1}}\right)\right]^{\frac{1}{p-2}}.
\end{equation}
It then follows from (\ref{clink2}) and 
(\ref{deflittler1}) that 
\begin{equation}\label{LinkLem1Eqn0005}
	J(z) \geqslant \frac{r^2}{4}
	\left(1-\frac{\mu_2}{\lambda_{m+1}}\right),
		\quad\mbox{ for } z\in X^+, 
			\mbox{ with } \|z\|_{X} =r.
\end{equation}
Thus, setting 
\begin{equation}\label{defa1}
    a_1=\frac{1}{4}\left(1-\frac{\mu_2}{\lambda_{m+1}}\right)
    \left[\frac{1}{4C}\left(1-\frac{\mu_2}{\lambda_{m+1}}\right)\right]^{\frac{2}{p-2}},
\end{equation}
we see from (\ref{LinkLem1Eqn0005}) that 
\begin{equation}\label{defrho1}
J(z) \geqslant a_1,
\quad\mbox{ for all }z\in X^{+},
    \mbox{ with } \|z\|_{X} =r,
\end{equation}
which yields the conclusion of the lemma in 
(\ref{cond1link}).
\end{proof}

Next, we present a second result that is needed in
the proof of existence of a nontrivial solution of the system (\ref{mainsys}) through application of Theorem \ref{linktheo}.

\begin{lemma}
Assume that  ($H_1$), ($H_2$), ($H_3$) and $(AR)_\theta$ are satisfied, and that the 
functional $J\colon X\to\R$ is as defined in 
 (\ref{func0}), where $Q$ is defined in (\ref{ARq0}). 
Let $X=X^{-}\oplus X^{+}$ be the decomposition of the Hilbert space $X$ given in (\ref{decz1}), (\ref{DfnDecompXEqn05}) and (\ref{dec1}). 
Then, for the choice of $r$ in (\ref{deflittler1}) and of  $a_1$ in (\ref{defa1}), 
there exist $\rho > r$ and $\delta_1 >0 $ such that  
$\lambda_{m}<\lambda_{m+1}-\delta_1 
< \mu_1\leqslant\mu_2<\lambda_{m+1}$, for 
$m \geqslant 1$, where $\mu_1$ and $\mu_2$ are the
eigenvalues of the matrix $\textbf{M}$ given in 
(\ref{eq1}), and   
\begin{equation}
 a_o:=\max_{z\in M_0} J(z) < a_1,
\label{cond1link55}
\end{equation}
where $M_0$ is the set 
\begin{equation}\label{M0Dfn05}
	M_0= \{z_1+t z_2 \mid z_1\in X^{-},\|z\|_{X}=\rho \mbox{ and } t\geqslant 0 \mbox{ or } \|z\|_{X}\leqslant \rho \mbox{ and }t=0\},
\end{equation}
with $z_2\in span(\varphi_{m+1})\times  span(\varphi_{m+1})$,  
where $\varphi_{m+1}$ is an eigenfunction corresponding 
to the eigenvalue $\lambda_{m+1}$,
and $\|z_2\|_{X}=r$. 
\label{lemmalink2}
\end{lemma}
\begin{proof} [Proof:] We follow a line of reasoning similar to that in \cite[Lemma $3.5$]{RabSuWang}. 

Let $z=z_1+t z_2 \in M_0$, where  $z_1\in X^{-}$, $z_2\in X^{+}$ with $\|z_2\|_{X}=r$, $t\geqslant 0$, and $r$ is defined in (\ref{deflittler1}). 
Then, using the definition of the functional $J$ given  in (\ref{func0}), the assumption that $F(x,z)\geqslant 0$ for all $z\in X$ and  $x\in\Omega$ given in ($H_3$), and the inequality (\ref{pest342}), we get 
the estimate
\begin{equation}\label{lm431}
	J(z)\leqslant 
	\frac{1}{2}\|z_{1}\|_{X}^2 + \frac{t^2}{2}\|z_2\|_{X}^2 -\frac{\mu_1}{2}\left(\|z_{1}\|_{Y}^2 + t^2\|z_{2}\|_{Y}^2\right).
\end{equation}
Next, write $z_1=(u_1,v_1)$ with $u_1,v_1\in X_m^{-}$. It then follows from the definition of $X_{m}^{-}$ 
in (\ref{dec1}) that 
\begin{equation}
    \|u_{1}\|^2\leqslant \lambda_m\|u_1\|_{L^2}^2\quad\mbox{ and }\quad \|v_1\|^2\leqslant \lambda_{m}\|v_1\|_{L^2}^2.
    \label{lm432}
\end{equation}
Substituting the estimates in  (\ref{lm432}) into (\ref{lm431}), we obtain
\begin{equation}
    J(z_1 + t z_2) \leqslant \left(\frac{\lambda_m}{2}-\frac{\mu_1}{2}\right)\|z_1\|_{Y}^2 + \frac{t^2}{2}\|z_2\|_{X}^2 -\frac{\mu_1 t^2}{2}\|z_2\|_{Y}^2.
    \label{lm433}
\end{equation}
Write $z_2=(s \varphi_{m+1},s\varphi_{m+1})$, for some $s\ne 0$. Then, using the definition of the inner product  in (\ref{inner3}), we obtain that 
\begin{equation}
    \|z_{2}\|_{Y}^2 = 2s^2\|\varphi_{m+1}\|_{L^2}^2.
    \label{lm434}
\end{equation}
Since $\|\varphi_{m+1}\|^2=\lambda_{m+1}\|\varphi_{m+1}\|_{L^2}^2$ and $\|\varphi_{m+1}\|=1$, it follows from (\ref{lm434}) that
\begin{equation}
    \|z_{2}\|_{Y}^2 = \frac{2}{\lambda_{m+1}}s^2.
    \label{lm435}
\end{equation}
On the other hand, since $\|z_{2}\|_{X}=r$, where $r$ is given by (\ref{deflittler1}),  
\begin{equation}
    r^2=\|z_{2}\|_{X}^2 = 2s^2\|\varphi_{m+1}\|^2=2s^2,
    \label{lm4391}
\end{equation}
for $s\ne 0$, since $\|\varphi_{m+1}\|=1.$ 

Combining (\ref{lm435}) and )\ref{lm4391}) then yields
\begin{equation}\label{lm43Eqn0005}
    \|z_{2}\|_{Y}^2 = \frac{r^2}{\lambda_{m+1}}.
\end{equation}
Hence, substituting (\ref{lm43Eqn0005}) into the estimate in (\ref{lm433}), and using the assumption that 
$\|z_2\|_{X} = r$, yields
\begin{equation}\label{lm436}
	J(z) \leqslant \frac{1}{2}(\lambda_m-\mu_1)\|z_1\|_{Y}^2 
		+ \frac{1}{2}
\left(1 - \frac{\mu_1}{\lambda_{m+1}}\right)
t^2r^2
\end{equation}
Given the assumption $\mu_1 > \lambda_{m}$, it follows from (\ref{lm436}) that 
\begin{equation}\label{lm437}
  J(z) \leqslant \frac{1}{2}
\left(1 - \frac{\mu_1}{\lambda_{m+1}}\right)
t^2r^2.
\end{equation}

Since we are assuming that $\lambda_{m+1}-\delta_1<\mu_1$, it follows from (\ref{lm437}) that 
\begin{eqnarray}\label{lm441}
    J(z) \leqslant \frac{\delta_1}{2\lambda_{m+1}} t^2r^2.
\end{eqnarray}

Next, set $\|z\|_{X}=\rho$ and use the definition of the norm of $X$ in (\ref{norm1}) to get
$$
    \rho^2 = \|z\|_{X}^2 = \|z_1\|_{X}^2+t^2r^2,
$$
from which we get that 
\begin{equation}\label{lm438}
     t^2r^2 \leqslant \rho^2.
\end{equation}
Combining the estimates in (\ref{lm441}) and (\ref{lm438})
then yields
\begin{eqnarray}\label{lm43Eqn0010}
    J(z) \leqslant \frac{\delta_1}{2\lambda_{m+1}} \rho^2,
\end{eqnarray}
for $z=z_1 + t z_2$, where $t\geqslant 0$, $z_1\in X^{-}$,
and $z_2=(s \varphi_{m+1},s\varphi_{m+1})$, for some $s\ne 0$, with $\|z_2\|_X = r$, and $\|z\|_X = \rho$.

Put 
\begin{eqnarray}\label{lm43Eqn0015}
    \delta_1 = \frac{ a_1\lambda_{m+1}}{\rho^2} ,
\end{eqnarray}
where $a_1$ is given by (\ref{defa1}).
Then, choosing $\rho > r$ sufficiently large so that 
$$
    \delta_1 < \lambda_{m+1} - \lambda_m,
$$
according to (\ref{lm43Eqn0015}), we obtain from 
(\ref{lm43Eqn0010}) and (\ref{lm43Eqn0015}) that 
\begin{eqnarray}\label{lm43Eqn0020}
    J(z) \leqslant \frac{a_1}{2},
\end{eqnarray}
for $z=z_1 + t z_2$, where $t\geqslant 0$, $z_1\in X^{-}$,
and $z_2=(s \varphi_{m+1},s\varphi_{m+1})$, for some $s\ne 0$, with $\|z_2\|_X = r$, and $\|z\|_X = \rho$.

On the other hand, 
for $z=z_1 + t z_2$, where $t = 0$, $z_1\in X^{-}$,
and $z_2=(s \varphi_{m+1},s\varphi_{m+1})$, for some $s\ne 0$, with $\|z_2\|_X = r$, and 
$\|z\|_X \leqslant \rho$, we obtain from (\ref{lm433}) that 
\begin{eqnarray}\label{lm43Eqn0025}
    J(z) \leqslant 0,
\end{eqnarray}
since $\lambda_m < \mu_1$.

Finally, combining (\ref{lm43Eqn0020}) and 
(\ref{lm43Eqn0025})  we obtain the assertion in 
(\ref{cond1link55}) in view of the definition of 
of $M_0$ in (\ref{M0Dfn05}). This concludes the proof of the lemma. 
\end{proof}

Next, we present the proof of existence of a nontrivial solution of the system (\ref{mainsys}).

 \begin{theorem}   
 Assume that  ($H_1$), ($H_2$), ($H_3$) and $(AR)_\theta$ are satisfied, and that the 
functional $J\colon X\to\R$ is as defined in 
 (\ref{func0}), where $Q$ is defined in (\ref{ARq0}). 
 For the Hilbert space  $X=H_{0}^{1}(\Omega)\times H_{0}^{1}(\Omega)$, consider 
 the decomposition $X=X^{-}\oplus X^{+}$, with 
 $$X^{-}=X_{m}^{-}\times X_{m}^{-}\quad\mbox{ and }\quad  X^{+}=(X^{-})^{\perp},$$ 
 where $X_m^{-}$ is the finite-dimensional space given by 
 \begin{equation}
    X_{m}^{-}=\bigoplus_{j\leqslant m}\ker(-\Delta - \lambda_{j}I).
    \nonumber 
\end{equation}
Suppose that there exists $m\geqslant 1$ such 
that 
$\lambda_m<\mu_1\leqslant\mu_2<\lambda_{m+1}$,
where $\mu_1$ and 
$\mu_2$ are the eigenvalues of the matrix $\textbf{M}$ given in
(\ref{eq1}).

Then, there exists $\delta_1>0$  such that $$\lambda_{m}<\lambda_{m+1}-\delta_1 < \mu_1\leqslant\mu_2<\lambda_{m+1}$$ 
and  the functional $J$ defined in (\ref{func0}) has at least one nontrivial critical point $z_o$.  
\label{maintheo4}
\end{theorem}

\begin{proof}[Proof:] By virtue of Lemma \ref{lemmalink1} and Lemma \ref{lemmalink2}, the functional $J$ satisfies the condition (\ref{bcond1}). Moreover, according to Lemma \ref{pslemma}, the functional $J$ also satisfies the Palais-Smale condition. 
 Therefore, it follows from the Linking Theorem \ref{linktheo} that $J$ has a nontrivial critical point $z_0$ corresponding to a critical value $c_{0}$ with
 $c_0\geqslant a_1 > 0$, where $a_1$  is given in Lemma \ref{lemmalink1}.   
\end{proof}

\section{Second Multiplicity result}\label{maintheosec}

In this section, we present a second multiplicity result for the system (\ref{mainsys}). We 
obtain a second nontrivial solution of the 
problem using an argument by contradiction 
involving the exact sequence of a specific topological pair. To get the second nontrivial critical point for $J$, in addition to the hypotheses $(H_1)$--$(H_3)$ and $(AR)_\theta$, we need to make the following assumption on the functions $f$ and $g$:
\begin{itemize}
    \item [$(H_4)$] $\displaystyle \frac{\partial f}{\partial s}(x,s,t),\displaystyle \frac{\partial f}{\partial t}(x,s,t),\displaystyle \frac{\partial g}{\partial s}(x,s,t),$ and $\displaystyle \frac{\partial g}{\partial t}(x,s,t)$ exist and are continuous on $\overline{\Omega}\times\R^2.$
\end{itemize}
In this case, $J\in C^{2}(X,\R)$. This will allow us to analyze the cases where the critical point $z_o$ obtained in  Theorem \ref{maintheo4} is a nondegenerate or degenerate critical point for $J$ by using the generalized Morse Lemma, also known as the splitting lemma (see \cite[Theorem $8.3$]{MW1}), and the shifting theorem (see \cite[Theorem $8.4$]{MW1}). 

First, we will need to compute the critical groups of $J$ at the origin and at infinity. Then, in the third subsection, we present the final argument.

\subsection{Critical groups of J at the origin} 

In this subsection, we will compute the critical groups of the functional $J$ at the origin. To do so, we will use a result presented in \cite{PerSch}, whose proof can be found in Chang and Ghoussub in \cite{ChGh}:
\begin{theorem}
\cite[Theorem $1.4.4$]{PerSch} Let $J_{t},$ $t\in[0,1],$ be a family of $C^{1}$ functionals on a Banach space $X$ and $z_{0}$ be a critical point of each $J_{t}$ within a closed neighborhood $U$.  Assume that
\begin{itemize}
    \item [(i)] $J_{t}$ satisfies the Palais-Smale condition over $U$, for each $t\in[0,1]$;
    \item [(ii)] $U$ contains no other critical points of $J_{t}$; and 
    \item [(iii)] the map $[0,1]\to C^{1}(U,\R)$: $t\mapsto J_{t}$ is continuous.
\end{itemize}
Then,
\begin{equation}
C_{q}(J_{0},z_{0})\cong C_{q}(J_{1},z_{0}),\quad\mbox{ for all }q\in\Z. 
\nonumber 
\end{equation}
\label{theoper}
\end{theorem}

\begin{remark}
{\rm Theorem (\ref{theoper}) implies that if the critical groups at an isolated critical point, $z_0$, for a particular functional $J_0$ are known, then  we can compute the critical groups at $z_0$ for a functional $J_{1}$ that can be continuously deformed into $J_0$ while preserving the Palais-Smale condition in an isolating neighborhood of $z_o$ for all functionals in the deformation.
}
\end{remark}

For $0\leqslant t\leqslant 1$, define $J_{t}:X\rightarrow\R$ to be  the family of $C^{1}$ functionals given by
\begin{equation}
    J_{t}(z)=Q(z)-t\int_{\Omega}F(x,z)\,dx,\quad\mbox{ for } z \in X,\mbox{ }t\in[0,1], 
    \label{jtdef}
\end{equation}
where $Q$ is the auxiliary functional discussed in Section \ref{auxFunc} and defined in (\ref{q2}).

The Fr\'{e}chet derivatives of the family of functionals defined in (\ref{jtdef}) are 
given by 
\begin{equation}
\langle J_{t}^{\prime}(z),\Lambda\rangle = \langle Q^{\prime}(z),\Lambda\rangle_{X} - t\int_{\Omega}\nabla F(x,z)\cdot \Lambda\,dx,
\label{jtfrechet}
\end{equation}
for all $\Lambda\in X$ and all $t\in[0,1]$. 

If we can show that the family of functionals $J_{t}$ defined in (\ref{jtdef}) satisfies the hypotheses of Theorem \ref{theoper}, then the critical groups of $J$ at the origin will be isomorphic to the critical groups of $Q$ at the origin as in (\ref{cgroupsQ}). 

We will need some auxiliary estimates. 

\begin{lemma}\label{gradestlemma0}
Let $F\colon\Omega\times\mathbb{R}^2\to\mathbb{R}$ be as given
by (\ref{gradFDef0}) and (\ref{gradFDef}), and assume that the 
hypotheses in ($H_1$) and ($H_2$) hold true. 
    \begin{itemize}
   \item [(i)] There exists a positive 
    constant $C$ such that 
\begin{equation}\label{8924Eqn0005}
    \left|\int_{\Omega}F(x,z)\,dx\right| 
             \leqslant C \|z\|_{X}^{p},
             	\quad\mbox{ for all } z\in X.
\end{equation}

\item [(ii)] There exists a positive 
    constant $C$ such that 
\begin{equation}\label{8924Eqn0010}
    \left|\int_{\Omega}\nabla F(x,z)\cdot \Phi\,dx\right| 
             \leqslant C \|z\|_{X}^{p-1}\|\Phi\|_{X},
             \quad\mbox{ for all } z,\Phi\in X.
\end{equation}
\end{itemize}
\end{lemma}

\begin{proof}[Proof:] 
    To establish the estimate (\ref{8924Eqn0005}) in assertion \textit{(i)},
    use (\ref{gradFDef0}) and (\ref{gradFDef})
    to write
\begin{equation}\label{AuxEst0005}
F(x,s,t) = 
\int_0^1 (f(x,rs, rt) s + g(x,rs,rt) t)\ dr,
    \quad\mbox{ for } x\in\Omega, (s,t)\in\R^2,
\end{equation}
where we have also used the chain rule.

Take absolute value on both sides of (\ref{AuxEst0005}) and 
apply the triangle inequality to obtain
\begin{equation}\label{AuxEst0010}
|F(x,s,t)| \leqslant 
\int_0^1 |f(x,rs, rt)| |s| \ dr  
	+ \int_0^1 |g(x,rs,rt)| |t|\ dr,
\end{equation}
for $x\in\Omega$ and $(s,t)\in\R^2$.

Next, use the assumption in ($H_2$) to obtain 
the estimate
$$
	\int_0^1 |f(x,rs, rt)| |s| \ dr 
		\leqslant
	C_1 \int_0^1 r^{p-1} (|s|^{p-1} + |t|^{p-1}) |s| \ dr,
$$
for $x\in\Omega$ and $(s,t)\in\R^2$, which 
integrates to 
\begin{equation}\label{AuxEst0015}
\int_0^1 |f(x,rs, rt)| |s| \ dr \leqslant 
\frac{C_1}{p}(|s|^p + |s||t|^{p-1}),
	\quad\mbox{ for } x\in\Omega \mbox{ and } (s,t)\in\R^2.
\end{equation}
Applying Young's inequality, 
$$
	ab \leqslant \frac{1}{p} a^p + \frac{1}{q} b^q,
		\quad\mbox{ for } a, b\geqslant 0, 
			p,q >1, 
			\frac{1}{p}+\frac{1}{q} =1,
$$
to the right-most term in (\ref{AuxEst0015}),
with $a=|s|$ and $b=|t|$, we obtain 
$$
\int_0^1 |f(x,rs, rt)| |s| \ dr \leqslant 
\frac{C_1}{p}(|s|^p + \frac{1}{p}|s|^p + \frac{p-1}{p}|t|^p),
$$
or
$$
\int_0^1 |f(x,rs, rt)| |s| \ dr \leqslant 
\frac{C_1}{p}(\frac{p+1}{p}|s|^p + \frac{p-1}{p}|t|^p),
$$
from which we get 
\begin{equation}\label{AuxEst0020}
\int_0^1 |f(x,rs, rt)| |s| \ dr \leqslant 
\frac{C_1(p+1)}{p^2}(|s|^p + |t|^{p}),
\end{equation}
for $x\in\Omega$ and $(s,t)\in\R^2$.

Similar calculations to those leading to 
(\ref{AuxEst0020}) can be used to obtain
\begin{equation}\label{AuxEst0025}
\int_0^1 |g(x,rs, rt)| |t| \ dr \leqslant 
\frac{C_2(p+1)}{p^2}(|s|^p + |t|^{p}),
\end{equation}
for $x\in\Omega$ and $(s,t)\in\R^2$.

Combining (\ref{AuxEst0010}), (\ref{AuxEst0020}) and 
(\ref{AuxEst0025}) then yields 
$$
|F(x,s,t)|  \leqslant 
\frac{(p+1)(C_1+C_2)}{p^2}(|s|^p + |t|^{p}),
$$
for $x\in\Omega$ and $(s,t)\in\R^2$,
or
\begin{equation}\label{AuxEst0030}
|F(x,s,t)|  \leqslant 
C_p(|s|^p + |t|^{p}),
\quad\mbox{ for } x\in\Omega \mbox{ and } (s,t)\in\R^2,
\end{equation}
and some positive constant $C_p$ that depends 
on $p$.

It follows from the estimate in (\ref{AuxEst0030}) that 
\begin{equation}\label{AuxEst0035}
    \int_{\Omega}|F(x,u,v)|\,dx 
             \leqslant C_p( \|u\|_{L^p}^{p} + \|v\|_{L^p}^{p}),
             	\quad\mbox{ for } u,v\in H^1_0(\Omega).
\end{equation}
Now, since $p<2^\ast$, it follows from 
(\ref{AuxEst0035}) and the Sobolev embedding 
theorem that 
\begin{equation}\label{AuxEst0040}
    \int_{\Omega}|F(x,u,v)|\,dx 
             \leqslant C( \|u\|^{p} + \|v\|^{p}),
             	\quad\mbox{ for } u,v\in H^1_0(\Omega),
\end{equation}
where $C$ is a positive constant and $\|\cdot\|$ is the norm 
in $H^1_0(\Omega)$ given in (\ref{ARnorm0}).  From this point on
in this paper, the symbol $C$ will denote various positive
constants that might have different values.

It follows from (\ref{AuxEst0040}), and the 
definition of the norm $\|\cdot\|_X$ in $X$
given in (\ref{norm2}), that 
$$
    \int_{\Omega}|F(x,u,v)|\,dx 
             \leqslant C\|z\|_X^{p},
             	\quad\mbox{ for } z=(u,v)\in X,
$$
and some positive constant $C$, which yields the estimate in 
(\ref{8924Eqn0005}).  This completes the proof of the assertion in
item \textit{(i)} of the lemma.

To prove the assertion in item \textit{(ii)} of
the lemma, write $z = (u,v)$ and 
$\Phi = (\varphi,\psi)$, and compute
\begin{equation}\label{AuxEst0045}
\left|\int_{\Omega}\nabla F(x,z)\cdot \Phi\,dx\right|
		\leqslant 
		\int_{\Omega}|f(x,u,v)||\varphi|\,dx + \int_{\Omega}|g(x,u,v)||\psi|\,dx .
\end{equation}
Next, apply H\"{o}lder's inequality with exponents $p$ and $q=p/(p-1)$ to obtain from 
 (\ref{AuxEst0045}) that
 \begin{equation}\label{AuxEst0050}
\begin{array}{rcl}
\displaystyle
	\left|\int_{\Omega}\nabla F(x,z)\cdot z\,dx\right|
		&\leqslant& \displaystyle
\left(\int_{\Omega}|f(x,z)|^{q}\,dx\right)^{1/q}\|\varphi\|_{L^{p}}\\
				\\
& & \displaystyle\qquad 
			+ \left(\int_{\Omega}|g(x,z)|^{q}\,dx\right)^{1/q}\|\psi\|_{L^{p}}.
\end{array}	
\end{equation}
Now, since $p < 2^{*}$, we can apply the Sobolev embedding theorem
to obtain from (\ref{AuxEst0050}) that there exists a positive constant $C$ such that 
\begin{equation}\label{AuxEst0055}
\left|\int_{\Omega}\nabla F(x,z)\cdot z\,dx\right|
  \leqslant  C\left(\|f(\cdot,z)\|_{L^q}
  + \|g(\cdot,z)\|_{L^q}\right)
  \|\Phi\|_{X}. 
\end{equation}

On the other hand, using the inequality $(a+b)^{q}\leqslant 2^{q-1}(a^q+b^q),$ for $a,b>0$ and $q\geqslant 1$, and the subcritical growth condition for $f$ and $g$ in (\ref{subcg}), see also assumption ($H_2$), we obtain that  
\begin{equation}\label{AuxEst0060}
    |f(x,z)|^{q}\leqslant 2^{q-1}C_1\left(|u|^{(p-1)q} + |v|^{(p-1)q}\right),
        \quad\mbox{ for } x\in\Omega,
\end{equation}
and
\begin{equation}\label{AuxEst0065}
   |g(x,z)|^{q}\leqslant 2^{q-1}C_2\left(|u|^{(p-1)q} + |v|^{(p-1)q}\right),
        \quad\mbox{ for } x\in\Omega,
\end{equation}
for positive constants $C_1$ and $C_2$ given in ($H_2$).

It follows from (\ref{AuxEst0060}) that
\begin{equation}\label{AuxEst0070}
\|f(\cdot, z)\|_{L^{q}}^{q} \leqslant C\left[\int_{\Omega}|u|^{p}\,dx  +   \int_{\Omega}|v|^{p}\,dx \right],
\end{equation}
where $C$ is a positive constant. 

Since $p < 2^{*},$ we can invoke the Sobolev embedding theorem to get from (\ref{AuxEst0070}) that 
\begin{equation}\label{AuxEst0075}
    \|f(\cdot, z)\|_{L^{q}}^{q} 
    \leqslant C\left( \|u\|^{p}+\|v\|^{p}\right),        
\end{equation}
for some positive constant $C$, where 
$\|\cdot\|$ is the norm in $H^1_0(\Omega)$
define by (\ref{ARnorm0}). Then, using the 
definition of norm $\|\cdot\|_X$ given in
(\ref{norm2}), we obtain from 
(\ref{AuxEst0075}) that 
$$
    \|f(\cdot, z)\|_{L^{q}}^{q} 
    \leqslant C  \|z\|_X^{p},        
$$
for some positive constant $C$;
consequently, 
\begin{equation}\label{AuxEst0080}
    \|f(\cdot, z)\|_{L^{q}} 
    \leqslant C  \|z\|_X^{p-1},        
\end{equation}
for some positive constant $C$, where we have also used the identity $\dfrac{p}{q} = p-1$.

Similar calculations to to those leading from 
(\ref{AuxEst0060}) to (\ref{AuxEst0080}) can be used to obtain from(\ref{AuxEst0065}) the estimate
\begin{equation}\label{AuxEst0085}
    \|g(\cdot,z)\|_{L^{q}} 
    \leqslant C\|z\|_{X}^{p-1},        
\end{equation}
for some positive constant $C$.

Substituting the estimates (\ref{AuxEst0080}) and (\ref{AuxEst0085}) into (\ref{AuxEst0055}), we obtain the estimate (\ref{8924Eqn0010}).  
This concludes the proof of the assertion in item \textit{(ii)} of the lemma. 
\end{proof}

\begin{lemma}\label{gradestlemma}
Let $X=X^{-}\oplus X^{+}$ be the decomposition of the Hilbert space $X$ defined in (\ref{decz1}), (\ref{DfnDecompXEqn05}) and
(\ref{dec1}). Let $F\colon\Omega\times\mathbb{R}^2\to\mathbb{R}$ by as given
by (\ref{gradFDef0}) and (\ref{gradFDef}), and assume that the 
hypotheses in ($H_1$) and ($H_2$) hold true.
For $z=z_{1}+z_2$ with $z_1\in X^{-}$ and $z_{2}\in X^{+}$, the following hold true.  
    \begin{itemize}
   \item [(i)] There exists a positive 
    constant $C$ such that 
\begin{equation}\label{32224Eqn0005}
    \left|\int_{\Omega}\nabla F(x,z)\cdot z\,dx\right| 
             \leqslant C \|z\|_{X}^{p}.
\end{equation}

\item [(ii)] There exists a positive 
    constant $C$ such that 
\begin{equation}\label{32224Eqn00051}
    \left|\int_{\Omega}\nabla F(x,z)\cdot z_{i}\,dx\right| 
             \leqslant C \|z\|_{X}^{p},
             \quad\mbox{ for }i=1,2.
\end{equation}

\end{itemize}
\end{lemma}

\begin{proof}[Proof:] 
 The assertion in item \textit{(i)}  follows
 from the assertion in item \textit{(ii)} 
 of Lemma \ref{gradestlemma0} by setting 
 $\Phi = z$.

To prove the assertion in item \textit{(ii)} of the lemma, first apply the assertion in 
item \textit{(ii)}  of Lemma \ref{gradestlemma0}
with $\Phi = z_i$, for $i=1,2$,  to obtain 
\begin{equation}\label{AuxEst0090}
\left|\int_{\Omega}\nabla F(x,z)\cdot z_i\,dx\right| 
\leqslant C \|z\|_{X}^{p-1}\|z_i\|_{X}, 
    \quad\mbox{ for } i=1,2,
\end{equation}
and some positive constant $C$.  Then, since
$\|z_i\|_{X}\leqslant\|z\|_{X}$, for $i=1,2$,
the assertion in item \textit{(ii)} of the
lemma follows from (\ref{AuxEst0090}).
This concludes the proof of the lemma.
\end{proof}

In what follows, we will show that the functional $J_{t}$, defined in (\ref{jtdef}), satisfies each of the  hypotheses of Theorem \ref{theoper}. 

We first show that there is an isolating neighborhood
of the origin for the family of functionals $J_t$ given in
(\ref{jtdef}) for all $t\in [0,1]$.

\begin{lemma}
Let $X=X^{-}\oplus X^{+}$ be the decomposition of the Hilbert space $X$ defined in (\ref{decz1}), (\ref{DfnDecompXEqn05}) and
(\ref{dec1}). Let $F\colon\Omega\times\mathbb{R}^2\to\mathbb{R}$ by as given by (\ref{gradFDef}) and (\ref{gradFDef0}), and assume that the 
hypotheses in ($H_1$) and ($H_2$) hold true.
 Let $J_t:X\rightarrow\R$ be the family of $C^1$ functionals defined in (\ref{jtdef}).  There exists
 a closed neighborhood, $U$, of the origin in $X$ such that the origin is the only critical point of $J_t$ in $U$, for all $t\in [0,1]$.
 \label{condii}
\end{lemma}

\begin{proof}[Proof:] 
Let $z=z_1+z_2$ with $z_1\in X^-$ and $z_2\in X^{+}$. Then, setting $\Phi=z_2-z_1$ in (\ref{jtfrechet}), we get 
 \begin{equation}
\langle J_{t}^{\prime}(z),z_2-z_1\rangle = \langle Q^{\prime}(z),z_2-z_1\rangle_{X} - t\int_{\Omega}\nabla F(x,z)\cdot (z_2-z_1)\,dx,
\label{jtfrechet1}
\end{equation}
for all $t\in[0,1]$. 

Substitute the estimate in (\ref{31024Eqn0020}) into (\ref{jtfrechet1}) to obtain
\begin{equation}
    \langle J_{t}^{\prime}(z),z_2-z_1\rangle  \geqslant \nu \|z\|_{X}^2- t\int_{\Omega}\nabla F(x,z)\cdot (z_2-z_1)\,dx,
\label{jtfrechet2}
\end{equation}
for all $t\in[0,1]$, where $\nu$ is given in
(\ref{31024Eqn0025}).

Applying the estimates in Lemma \ref{gradestlemma} we obtain the estimate
\begin{equation}
    \left|\int_{\Omega}\nabla F(x,z)\cdot (z_2-z_1)\,dx \right|\leqslant C\|z\|_{X}^{p},
    \label{gradest2}
\end{equation}
where $C$ is a positive constant.

Hence, it follows from (\ref{jtfrechet2}) and (\ref{gradest2}) that 
\begin{equation}
    \langle J_{t}^{\prime}(z),z_2-z_1\rangle  \geqslant \nu \|z\|_{X}^2\left(1-  \frac{C}{\nu}\|z\|_{X}^{p-2}\right),
    \label{jtfrechet3}
\end{equation}
for $z=z_1 + z_2 \in X = X^-\oplus X^+$, 
where $C$ is a positive constant and $\nu$ is given in (\ref{31024Eqn0025}).

Put 
\begin{equation}\label{ARisolEqn0005}
\rho=\left(\frac{\nu}{2C}\right)^{1/(p-2)}.
\end{equation}
Then, we obtain from (\ref{jtfrechet3}) that 
\begin{equation}\label{jtfrechet5}
    \langle J_{t}^{\prime}(z),z_2-z_1\rangle  \geqslant \frac{\nu}{2} \|z\|_{X}^2,
        \quad\mbox{ for } \|z\|_X \leqslant \rho.
\end{equation}
Hence, if $\overline{z}$ is a critical point of $J_t$, and $\|\overline{z}\| \leqslant \rho$, it follows from (\ref{jtfrechet5}) that
$$0\geqslant \frac{\nu}{2}\|\overline{z}\|_{X}^2,$$
which implies that $\overline{z}=0$.  

Thus, setting 
$$
    U = \overline{B_\rho(0)}
        = \{ z\in X\mid \|z\|\leqslant \rho\}
$$
where $\rho$ is given in (\ref{ARisolEqn0005}),
we see that $U$ is a closed neighborhood of the origin in $X$ such that the origin is an isolated critical point of $J_t$, for all $t\in [0,1].$
\end{proof}

 \begin{lemma}  
 Assume that the hypotheses in ($H_1$), ($H_2$), ($H_3$) and 
$(AR)_\theta$ hold true.
 The functional $J_{t}$ defined in (\ref{jtdef}) satisfies the Palais-Smale condition for all $t\in[0,1]$. 
 \label{pslemma1}
  \end{lemma} 

\begin{proof}[Proof:]
The case $t=0$ was proved in Lemma  \ref{psQ0lemma}.  The argument presented in  
the proof of Lemma \ref{pslemma} can be used to prove the case $t\in(0,1]$.
\end{proof}

Finally, we show that the family of functionals $J_t$ satisfy condition $(iii)$ of Theorem \ref{theoper}.

\begin{lemma}
Assume that the hypotheses in ($H_1$) and ($H_2$) hold true.
 Let $J_t:X\rightarrow\R$ be the family of $C^1$ functionals defined in (\ref{jtdef}).
Let $U = \overline{B_{\rho}(0)}$, where $\rho$ is given in (\ref{ARisolEqn0005}). Then, the function $\mathcal{F}:[0,1]\rightarrow C^{1}(U,\R)$ defined by 
\begin{equation}\label{FJtDfn}
	\mathcal{F}(t) = J_{t}
    \quad\mbox{ for } t\in [0,1],
\end{equation}
is continuous.  
\label{condiii}
\end{lemma}
\begin{proof}[Proof:] 
Endow $C^{1}(U,\R)$ with the norm 
$\|\cdot\|_{C^1(U)}$ given by 
\begin{equation}\label{condiiiEqn0005}
	\|J\|_{C^1(U)} = \sup_{z\in U} |J(z)| 
		+ \sup_{z\in U}\|J'(z)\|_{_{\mbox{op}}},
			\quad\mbox{ for } J\in C^1(U,\R),
\end{equation}
where $\|\cdot\|_{_{\mbox{op}}}$ is the operator norm in $X^\ast$, the dual space of $X$.

For $t_1$ and $t_2$ in $[0,1]$, use 
(\ref{condiiiEqn0005}) to compute 
\begin{equation}\label{condiiiEqn0010}
\|J_{t_2} - J_{t_1}\|_{C^1(U)} 
= \sup_{z\in U} |J_{t_2}(z) - J_{t_1}(z)| 
+ \sup_{z\in U}\|J_{t_2}'(z)
    -J_{t_1}'(z)\|_{_{\mbox{op}}},
\end{equation}
where, according to the definition of $J_t$ in
(\ref{jtdef}), 
$$
J_{t_2}(z) - J_{t_1}(z)=(t_1-t_2)\int_{\Omega}F(x,z)\,dx,\quad\mbox{ for } z \in X;  
$$
so that, 
\begin{equation}\label{condiiiEqn0015}
|J_{t_2}(z) - J_{t_1}(z) | 
\leqslant
|t_2-t_1| 
\left|\int_{\Omega}F(x,z)\,dx\right|,
    \quad\mbox{ for } z \in X.
\end{equation}

It follows from (\ref{condiiiEqn0015}) and the
estimate in (\ref{8924Eqn0005}) in Lemma 
\ref{gradestlemma0} that
$$
|J_{t_2}(z) - J_{t_1}(z) | 
\leqslant
C\rho^p |t_2-t_1|,
    \quad\mbox{ for } z \in U,
$$ 
for a positive constant $C$; so that,
\begin{equation}\label{condiiiEqn0020}
\sup_{z\in U} |J_{t_2}(z) - J_{t_1}(z) | 
\leqslant C\rho^p |t_2-t_1|,
\end{equation} 
for a positive constant $C$.

Similarly, using the definition of the derivative
of $J_t$ in (\ref{jtfrechet}), we obtain that 
$$
\langle J_{t_2}^{\prime}(z)-J_{t_1}^{\prime}(z),\Lambda\rangle 
= (t_1 - t_2) \int_{\Omega}\nabla F(x,z)\cdot \Lambda\,dx,
    \quad\mbox{ for } z , \Lambda \in X.
$$
Consequently,
\begin{equation}\label{condiiiEqn0025}
|\langle J_{t_2}^{\prime}(z)-J_{t_1}^{\prime}(z),\Lambda\rangle|
\leqslant 
|t_2 - t_1|
\left|\int_{\Omega}\nabla F(x,z)\cdot \Lambda\,dx\right|,
    \quad\mbox{ for } z , \Lambda \in X.
\end{equation}
It then follows from (\ref{condiiiEqn0025}) and
the estimate in (\ref{8924Eqn0010}) in the assertion in item \textit{(ii)} of
Lemma \ref{gradestlemma0} that
$$
|\langle J_{t_2}^{\prime}(z)-J_{t_1}^{\prime}(z),\Lambda\rangle|
\leqslant 
C \|z\|_X^{p-1} \|\Lambda\|_X |t_2 - t_1|,
    \quad\mbox{ for } z , \Lambda \in X,
$$
for some positive constant $C$; so that,
\begin{equation}\label{condiiiEqn0030}
|\langle J_{t_2}^{\prime}(z)-J_{t_1}^{\prime}(z),\Lambda\rangle|
\leqslant 
C \rho^{p-1} |t_2 - t_1|,
\end{equation}
for $z \in U$ and $\Lambda \in X$ with\
$\|\Lambda\|_X\leqslant 1$.

It then follows from (\ref{condiiiEqn0030}) that
$$
\sup_{\|\Lambda\|_X\leqslant 1} 
|\langle J_{t_2}^{\prime}(z)-J_{t_1}^{\prime}(z),\Lambda\rangle|
	\leqslant 
C \rho^{p-1} |t_2 - t_1|,
    \quad\mbox{ for } z \in U ;
$$
so that,
\begin{equation}\label{condiiiEqn0035}
\sup_{z\in U}\|J_{t_2}'(z)-J_{t_1}'(z)\|_{_{\mbox{op}}}
\leqslant 
C \rho^{p-1} |t_2 - t_1|.
\end{equation}

Combining (\ref{condiiiEqn0010}) with the estimates in 
(\ref{condiiiEqn0020}) and (\ref{condiiiEqn0035}), we obtain 
the estimate
\begin{equation}\label{condiiiEqn0040}
\|J_{t_2} - J_{t_1}\|_{C^1(U)} 
\leqslant 
C(\rho) |t_2 - t_1|,
\end{equation}
where $C(\rho)$ is a positive constant that depends on 
$\rho$.

The continuity of the map
$\mathcal{F}:[0,1]\rightarrow C^{1}(U,\R)$ 
given in (\ref{FJtDfn}) now follows from the 
estimate in (\ref{condiiiEqn0040}).
\end{proof} 

We are now ready prove the main result of this
section.
\begin{theorem}\label{CqJ0Thm}
Assume that there exists 
$m\geqslant 1$ such that 
$\lambda_m<\mu_1\leqslant\mu_2<\lambda_{m+1}$,
where $\mu_1$ and $\mu_2$ are the eigenvalues of the matrix 
$\textbf{M}$ given in (\ref{eq1}), and let $X=X^{-}\oplus X^{+}$ be the decomposition of the Hilbert space $X$ defined in (\ref{decz1}), (\ref{DfnDecompXEqn05}) and
(\ref{dec1}), where $X^{-}$ is finite dimensional. 
Let
$J\colon X\to\R$ is as defined in 
 (\ref{func0}), where $Q$ is defined in (\ref{ARq0}), and  
assume that the hypotheses in ($H_1$), ($H_2$), ($H_3$) and 
$(AR)_\theta$ hold true.
Then, $0$ is an isolated critical point of $J$ and the critical
groups of $J$ at $0$ are given by
\begin{equation}\label{CritGrpJ0}
    C_{k}(J,0)\cong \delta_{k,\dim X^{-}}\F,\quad\mbox{ for all }k\in\Z.
\end{equation}
\end{theorem} 

\begin{proof}[Proof:]
For $t\in[0,1]$, let $J_t \colon X\to\R$ be as
defined in (\ref{jtdef}).  By Lemma \ref{condii}, there exists a closed neighborhood
$U$ of the origin such that the origin is the only critical point of $J_t$ in $U$, for all $t\in [0,1].$ Thus, $U$ is an isolating 
neighborhood of origin for the functional $J$. 
It follows from Lemma \ref{condiii} that the map 
$[0,1]\to C^{1}(U,\R)$: $t\mapsto J_{t}$ is continuous. 
By Lemma \ref{pslemma1}, the functional $J_{t}$ defined in (\ref{jtdef}) satisfies the Palais-Smale condition for all $t\in[0,1]$. 
Hence, the functionals $J_t$, for $t\in[0,1]$, satisfy the hypotheses of Theorem \ref{theoper}. Therefore, it follows from  (\ref{cgroupsQ}) that 
\begin{equation}
    C_{k}(J,0)\cong C_{k}(Q,0)\cong\delta_{k,\dim X^{-}}\F,\quad\mbox{ for all }k\in\Z, 
    \label{cgrouporigin} 
\end{equation}
where  $X^{-}=X_{m}\times X_{m}$ as given 
in (\ref{dec1}).
\end{proof}

\subsection{Critical Groups of J at Infinity} 
In this section, we compute the critical groups of $J$ at infinity. These critical groups were first  defined by Bartsch and Li in  \cite{BLi}. 

Let $\mathcal{K}$ denote the set of critical 
points of $J$. 
For the case in which the set of critical values of $J\in C^1 (X,\R)$ is bounded from below and ${J}$ satisfies the Palais-Smale   condition, the global behavior of ${J}$ can be described by the critical groups of $J$ at infinity, 
which are defined  by
\begin{equation}
    C_{q}(J,\infty) = H_{q}(X,J^{a_{0}}),\quad\mbox{ for all }q\in\Z,
    \label{cinfty}
\end{equation}
where $a_{0} < \inf J(\mathcal{K}).$ These groups are well--defined as a consequence of the second deformation lemma (see \cite[Lemma $1.1.2$]{PerSch}).

First, we show that that the set of critical points, $\mathcal{K}$, of $J$ defined in (\ref{func0}), where $Q$ is defined in (\ref{q2}), is bounded from below. 

\begin{lemma}
Let $J\colon X\to\R$ be as defined in 
(\ref{func0}), where $Q$ is defined in (\ref{q2}).  Assume that the hypotheses in ($H_1$), ($H_2$) and 
$(AR)_\theta$ hold true. Then, the set of critical values of $J$ is bounded from below.  \label{bbelow}
\end{lemma} 

\begin{proof}[Proof:]
Let  $z_o \in\mathcal{K}$, and 
Put $z=z_o$ and $\Lambda=z_o$ in (\ref{DerJ05}), 
to obtain 
 \begin{equation}
     \|z_o\|_{X}^2 - \int_{\Omega}\textbf{M}z_o\cdot z_o\,dx = \int_{\Omega}\nabla F(x,z_o)\cdot z_o\,dx,
     \label{secinf1} 
 \end{equation}
 since $\langle J'(z_o),z_o\rangle = 0$.

Next, use the definition of $J$ in (\ref{func0}) 
to obtain from (\ref{secinf1}) that 
\begin{equation}\label{CritGrpInf0005}
    J(z_o) = \int_{\Omega}\left(\frac{1}{2}\nabla F(x,z_o)\cdot z_o 
        - F(x, z_o)\right)\,dx.
\end{equation}

Let $R$ be as given in hypothesis $(AR)_{\theta}$ and put
\begin{equation}\label{CritGrpInf0010}
	\Omega_R = 
	\{ x\in\Omega\mid |z_o(x)|\leqslant R\}, 
\end{equation}
where $|\cdot|$ denotes the Euclidean norm in $\R^2$. It then
follows from the assumption in $(AR)_{\theta}$ that 
\begin{equation}\label{CritGrpInf0015}
	0 < \theta F(x,z_o(x)) < \grad F(x,z_o(x))\cdot z_o(x),
		\quad\mbox{ for } x\in\Omega_R^c, 
\end{equation}
where $\Omega_R^c$ denotes the complement of $\Omega_R$.

Use the definition of $\Omega_R$ in (\ref{CritGrpInf0010})
to rewrite (\ref{CritGrpInf0005}) as follows:
\begin{equation}\label{CritGrpInf0020}
\begin{array}{rcl}
	 J(z_o) & = & \displaystyle
	 \int_{\Omega_R}\left(\frac{1}{2}\nabla F(x,z_o)\cdot z_o 
        - F(x, z_o)\right)\,dx\\
        					\\
        	& &\displaystyle \qquad + \ 
        	\int_{\Omega_R^c}\left(\frac{1}{2}\nabla F(x,z_o)\cdot z_o 
        - F(x, z_o)\right)\,dx.
\end{array}
\end{equation}
Put
$$
C_1 = 
\max_{\genfrac{}{}{0pt}{2}{x\in\overline{\Omega}}{|(s,t)|\leqslant R}}
\left|\frac{1}{2} \grad F(x,s,t)\cdot (s,t) - F(x,s,t)\right|.
$$
It then follows from (\ref{CritGrpInf0020}) that 
\begin{equation}\label{CritGrpInf0025}
J(z_o) \geqslant - C_1 |\Omega| +
	\int_{\Omega_R^c}\left(\frac{1}{2}\nabla F(x,z_o)\cdot z_o 
        - F(x, z_o)\right)\,dx,
\end{equation}
where $|\Omega|$ denotes the Lebesgue measure of $\Omega$.

The estimate in (\ref{CritGrpInf0015}) and the
assumption that $\theta > 2$ can be
used to show that 
$$
\frac{1}{2}\nabla F(\cdot,z_o)\cdot z_o - F(\cdot, z_o)
	> \frac{1}{\theta}\nabla F(\cdot,z_o)\cdot z_o - F(\cdot, z_o) > 0
		\quad\mbox{ in } \Omega_R^c.
$$
Therefore, we obtain from (\ref{CritGrpInf0025}) the estimate
$$
J(z_o)\geqslant  -C_1|\Omega|,
	\quad\mbox{ for all } z_o\in\mathcal{K},
$$
which shows that the set of critical values
of $J$ is bounded from below. This concludes the proof of the lemma. 
\end{proof} 

In what follows, let $a_{0}>0$ be such that
\begin{equation}
    -a_{0}<\displaystyle\inf_{z\in \mathcal{K}}J(z).
    \nonumber 
\end{equation} 
 
By following an argument similar to in Wang \cite[Section $3$]{Wang}, we will compute the critical groups of $J$ at infinity. First, we prove the following auxiliary results.

\begin{lemma}
Assume that the hypotheses in ($H_1$), ($H_2$) and $(AR)_\theta$ hold true, and let
$J\colon X\to\R$ be as defined in (\ref{func0}), 
where $Q$ is defined in (\ref{q2}). 
Then, for $z\in X\backslash\{0\}$, we have
\begin{equation}
    \lim_{t\rightarrow +\infty} J(tz)=-\infty.
    \label{jinf}
\end{equation}
\label{cinfl1}

\end{lemma}
\begin{proof}[Proof:]
Let $z\in X\backslash\{0\}$ and use (\ref{func0}) to
compute 
\begin{equation}\label{CritGrpInf0030}
	J(tz) = 
	Q(tz) - \int_{\Omega}F(x,tz)\,dx,
			\quad\mbox{ for } t > 0.
\end{equation}
Next, use the fact that $Q$ is a quadratic form,
and the estimate in (\ref{arr2}) to obtain 
from (\ref{CritGrpInf0030}) that 
\begin{equation}\label{CritGrpInf0035}
J(tz) \leqslant t^2Q(z)
-c_1 t^{\theta}\int_{\Omega}|z|^{\theta}\,dx 
  + c_2|\Omega|, 
        \quad\mbox{ for } t > 0,
\end{equation}
where $c_1$ and $c_2$ are positive constants. 
The assertion in (\ref{jinf}) now follows
from (\ref{CritGrpInf0035}) and the assumption
that $\theta > 2$ because $z\not=0$.
\end{proof} 

\begin{lemma}\label{CritGrpInfL3}
Assume that the hypotheses in ($H_1$), ($H_2$) and $(AR)_\theta$ hold true, and let
$J\colon X\to\R$ be as defined in (\ref{func0}), 
where $Q$ is defined in (\ref{q2}). 
Then, there exists $M_o>0$ such that, for all $M\geqslant M_o$,
\begin{equation}\label{CritGrpInfL3conclusion}
	\mbox{ if } (t,z)\in(0,\infty)\times (X\backslash\{0\})
	\mbox{ and } J(tz)\leqslant -M, \mbox{ then }
	\frac{d}{dt}[J(tz)] < 0.
\end{equation}
\end{lemma}

\begin{proof}[Proof:]
Let $z\in X\backslash\{0\}$ and $t>0$, and use 
(\ref{func0}) to compute 
$$
	J(tz) = 
	Q(tz) - \int_{\Omega}F(x,tz)\,dx;
$$
so that,  using the fact that $Q$ is a quadratic form, 
\begin{equation}\label{CritGrpInf0040}
	J(tz) = 
	t^2Q(z) - \int_{\Omega}F(x,tz)\,dx,
		\quad\mbox{ for } t >0  \mbox{ and } z\in X\backslash\{0\}.
\end{equation}
Next, use the chain rule to obtain from (\ref{CritGrpInf0040})
that
$$
	\frac{d}{dt} [J(tz)] = 
	2tQ(z)- \int_{\Omega}\nabla F(x,tz)\cdot z\,dx, 
	\quad\mbox{ for } t > 0 \mbox{ and } z\in X\backslash\{0\},
$$
which we can rewrite as
\begin{equation}\label{CritGrpInf0045}
	\frac{d}{dt} [J(tz)] = 
	\frac{2}{t}\left[t^2Q(z) - \frac{1}{2}\int_{\Omega} \nabla F(x,tz)\cdot (tz)\,dx   \right],
\end{equation}
for $t > 0$ and $z\in X\backslash\{0\}$.

Next, substitute (\ref{CritGrpInf0040}) into the 
right-hand side of (\ref{CritGrpInf0045}) to 
obtain
\begin{equation}\label{CritGrpInf0050}
	\frac{d}{dt} [J(tz)] = 
	\frac{2}{t}\left[J(tz) - 
	\int_{\Omega}\left(\frac{1}{2} \nabla F(x,tz)\cdot (tz)
	- F(x,tz)\right)\,dx   \right],
\end{equation}
for $t > 0$ and $z\in X\backslash\{0\}$.

For $t > 0$ and $z\in X\backslash\{0\}$, put
\begin{equation}\label{CritGrpInf0055}
	\Omega_{R,t} = 
	\{ x\in\Omega\mid |tz(x)|\leqslant R\}. 
\end{equation}
Then, using the assumption in $(AR)_{\theta}$, we get that 
\begin{equation}\label{CritGrpInf0060}
	0 < \theta F(x,tz(x)) < \grad F(x,tz(x))\cdot(t z(x)),
		\quad\mbox{ for } x\in\Omega_{R,t}^c, 
\end{equation}
where $\Omega_{R,t}^c$ denotes the complement of $\Omega_{R,t}$
defined in (\ref{CritGrpInf0055}).

Put
\begin{equation}\label{CritGrpInf0065}
C_o = 
\max_{\genfrac{}{}{0pt}{2}{x\in\overline{\Omega}}{w\in\R^2, |w|\leqslant R}}
\left|\frac{1}{2} \grad F(x,w)\cdot w - F(x,w)\right|.
\end{equation}
Then, assuming that $t>0$, 
$z\in X\backslash\{0\}$, and
$J(tz)\leqslant -M$, we obtain from 
(\ref{CritGrpInf0050}) that 
\begin{equation}\label{CritGrpInf0070}
\begin{array}{rcl}
	\displaystyle
	\frac{d}{dt} [J(tz)] 
	& \leqslant & \displaystyle
	\frac{2}{t}(-M + C_o|\Omega|) \\
				\\
	&  & \displaystyle \qquad
		-\ \frac{2}{t}
	\int_{\Omega_{R,t}^c}\left(\frac{1}{2} \nabla F(x,tz)\cdot (tz)
	- F(x,tz)\right)\,dx,   
\end{array}
\end{equation}
where $C_o$ is the constant given in (\ref{CritGrpInf0065}).

Let $M_o = 2C_o|\Omega|$  and suppose that 
$M\geqslant M_o$.  Then, it follows from 
(\ref{CritGrpInf0070}), the assertion in 
(\ref{CritGrpInf0060}) and the assumption that
$\theta > 0$, that 
$$
\frac{d}{dt} [J(tz)] 
\leqslant
-\frac{C_o|\Omega|}{t} <0.
$$
This establishes the assertion in 
(\ref{CritGrpInfL3conclusion}) and the proof
of the lemma is now complete.
\end{proof} 

 The next lemma will allow us to compute the singular homology groups of the level sets $J^{-M}=\{z\in X:J(z)\leqslant -M\}$ for $M\geqslant M_o$, where $M_o$ is given by Lemma 
 \ref{CritGrpInfL3}.

\begin{lemma}
Assume that the hypotheses in ($H_1$), ($H_2$) and $(AR)_\theta$ hold true.  Assume also that
there exists $m\geqslant 1$ such that $\lambda_m < \mu_1\leqslant\mu_2 < \lambda_{m+1}$. Let
$J\colon X\to\R$ be as defined in (\ref{func0}), 
where $Q$ is defined in (\ref{q2}). Let $M_o$ be 
as given by Lemma \ref{CritGrpInfL3}.
Then, there exists $M_1\geqslant M_o$ such that, for all $M\geqslant M_o$, $J^{-M}$ is homotopically equivalent to $S^{\infty}$, 
where $S^{\infty}=\{z\in X:\|z\|_{X}=1\}$ is the unit sphere in $X$.
\label{cinflema}
\end{lemma}

\begin{proof}[Proof:]
Let $M_o$ be as given by Lemma 
\ref{CritGrpInfL3} and let $M\geqslant M_o$.
For $z\in X\backslash\{0\}$, since $J(0) = 0$, it follows from the result of Lemma \ref{cinfl1} in
(\ref{jinf}) and the intermediate value theorem
that there exists $t>0$ such that 
\begin{equation}\label{CritGrpInf0075}
    J(tz) = -M, 
        \quad\mbox{ for } z\in S^{\infty}.
\end{equation}
It follows from (\ref{CritGrpInf0075}) and the 
result of Lemma \ref{CritGrpInfL3} in 
(\ref{CritGrpInfL3conclusion}), for the value 
of $t>0$ in (\ref{CritGrpInf0075}),
\begin{equation}\label{CritGrpInf0080}
    \frac{d}{dt}[J(tz)] < 0,
        \quad\mbox{ for } z\in S^{\infty}.
\end{equation}
In view of (\ref{CritGrpInf0080}) and the 
implicit function theorem (see \cite[Theorem $15.1$]{Deim1}), $t = T(z)$, where 
$T\in C^1(S^{\infty},\R)$; so that,  
\begin{equation}\label{CritGrpInf0085}
\begin{aligned}
    J(tz) & > -M \quad \mbox{ for } t< T(z);\\
    J(tz) & = -M \quad \mbox{ for }t=T(z),\\ 
    J(tz) & < -M \quad \mbox { for }t>T(z).
\end{aligned}   
\end{equation}

We show that, for $M$ large enough, we may 
assume that 
\begin{equation}
    T(z)\geqslant 1,\quad\mbox{ for all }z\in S^{\infty}.
    \label{lowerboundT}
\end{equation}
Indeed, using the definition of the functional $J$ in (\ref{func0}) we have
\begin{equation}
    \int_{\Omega}F(x,tz)\,dx = t^2Q(z)-J(tz), \quad\mbox{ for }z\in X; 
    \nonumber 
\end{equation}
so that, for $z\in S^{\infty}$ and $t=T(z)$, it follows from (\ref{8924Eqn0005}) and (\ref{CritGrpInf0085}) that
\begin{equation}
    C (T(z))^{p}\geqslant (T(z))^{2}Q(z) + M,\quad\mbox{ for }z\in S^{\infty}.
    \label{jtzineq2}
\end{equation}


Let $z=z_1+z_2\in X^{-}\oplus X^{+}$. It follows from the definition of $Q$ in (\ref{q2})  that
\begin{equation}
    Q(z) = \left[\frac{1}{2}\|z_1\|_{X}^2 - \frac{1}{2}\int_{\Omega}\textbf{M}z_1\cdot z_1\,dx\right] + \left[\frac{1}{2}\|z_2\|_{X}^2-\frac{1}{2}\int_{\Omega}\textbf{M}z_2\cdot z_2\,dx\right] .
    \label{jtz22}
\end{equation}

By virtue of (\ref{qd12}), we obtain from (\ref{jtz22}) that  
\begin{equation}
\begin{aligned}
    Q(z) & \geqslant \left[\frac{1}{2}\|z_1\|_{X}^2 - \frac{1}{2}\int_{\Omega}\textbf{M}z_1\cdot z_1\,dx\right] + \frac{1}{2}\left[1-\frac{\mu_2}{\lambda_{m+1}}\right]\|z_{2}\|_{X}^2 \\
    & \geqslant   - \frac{\mu_2}{2}\|z_{1}\|_{Y}^2,
    \end{aligned}
    \label{jtz23}
\end{equation}
for all $z\in S^{\infty}$, since 
$\mu_2 < \lambda_{m+1}.$ 

Next, observe that the space $X^{-}$ is finite-dimensional. Hence, there exists $C_1>1$ such that 
\begin{equation}
    \frac{1}{C_1}\|z_1\|_{X}^2 \leqslant \|z_1\|_{Y}^2 \leqslant C_1\|z_1\|_{X}^2,\quad\mbox{ for }z_1\in X^{-}.
    \label{equivnorm33}
\end{equation}
Thus, substituting (\ref{equivnorm33}) into (\ref{jtz23}), we get 
$$
    Q(z) \geqslant  - \frac{C_1\mu_2}{2}\|z_1\|_{X}^2 ;
$$
so that, since  $\|z_1\|_{X}\leqslant \|z\|_{X}=1$, given that $z\in S^{\infty}$,
\begin{equation}\label{jtz25}
    Q(z) \geqslant - C_2,
\end{equation}
where we have set $C_2=\dfrac{C_1\mu_2}{2}$. 

It follows from the estimates in (\ref{jtz25}) and (\ref{jtzineq2}) that 
\begin{equation}
    (T(z))^{p}\geqslant \frac{M}{C} -\frac{C_2}{C}(T(z))^2   ,\quad\mbox{ for }z\in S^{\infty};
    \nonumber 
\end{equation}
or, 
\begin{equation}
    (T(z))^{p} + C_{3}(T(z))^2 \geqslant \frac{M}{C},\quad\mbox{ for }z\in S^{\infty},
    \label{mainTeq} 
\end{equation}
where $C_3 = \dfrac{C_{2}}{C}$.

If 
\begin{equation}\label{CritGrpInf0090}
	0 < T(z) < 1,
		\quad\mbox{ for all } z\in S^\infty, 
\end{equation}
holds true, then  $(T(z))^{2} < T(z)$ and 
$(T(z))^{p} < T(z)$, since we are assuming that 
$p>1$.
It then follows from (\ref{mainTeq}) that 
$$
    (1+C_3)T(z) > \frac{M}{C};
$$
so that,
\begin{equation}
    T(z) > \frac{M}{C+C_2}, 
        \quad\mbox{ for }  z\in S^\infty.
    \label{test1}
\end{equation}
Observe from (\ref{test1}) that, if $M>C+C_2$,
then $T(z) >1$, which is in direct contradiction
with (\ref{CritGrpInf0090}). Consequently, 
if $M> C+ C_2$, the assertion in 
(\ref{lowerboundT}) must hold true. 

Let 
$$
	M_1 = \max\{M_o, C+C_2\}.
$$
Then, if $M\geqslant M_1$, 
\begin{equation}
	T(z)\geqslant 1,\quad\mbox{ for all }z\in S^{\infty}.
 \nonumber 
\end{equation}

Put 
$$
    B^\infty = \{ z\in X \mid \|z\|_Z <1\},
$$
the unit ball in $X$, and define 
$\eta:[0,1]\times(X\backslash B^\infty)\rightarrow X\backslash B^\infty$ as follows:
\begin{equation}
    \eta(s,z) = (1-s)z + sT(\pi_{S}(z))\pi_{S}(z),\quad\mbox{ for }z\in X\backslash B^\infty,
    \nonumber 
\end{equation}
where $\pi_{S}:X\backslash\{0\}\rightarrow S^{\infty}$ denotes the radial projection onto $S^{\infty}$; so that, 
\begin{equation}
    \pi_{S}(z)=\frac{1}{\|z\|_{X}}z;\quad\mbox{ for all }z\in X\backslash\{0\}.
    \nonumber 
\end{equation}
Then, $\eta$ is continuous. We also have that 
 $$\|\eta(s,z)\|_{X} =
    (1+s)\|z\|_{X} + s T(\pi_{S}(z))
 \geqslant 1,
    \quad\mbox{ for all }z\in X\backslash B^\infty.$$
Hence, $\eta$ does indeed map into $X\backslash B^\infty$. 

Note that $\eta(0,z)=z$ for all $z\in X\backslash B^\infty$ and $\eta(1,z)=T(\pi_{s}(z))\pi_{s}(z),$ for all $z\in X\backslash B^\infty.$ 
Therefore, we have 
\begin{equation}
    J(\eta(1,z))=-M,\quad\mbox{ for all }z\in X\backslash B^\infty.
    \nonumber 
\end{equation}
Consequently, $\eta(1,z)\in J^{-M}$ for all $z\in X\backslash B^\infty$. Thus, $J^{-M}$ is a deformation retract of $X\backslash B^\infty.$ Since $X\backslash B^\infty$ is homotopically equivalent to $S^{\infty}$, we conclude that 
\begin{equation*}
    J^{-M}\cong S^{\infty};
\end{equation*}
that is, $J^{-M}$ is homotopically equivalent to $S^{\infty}$ for $M\geqslant M_1$. This concludes the proof of the lemma. 
\end{proof}

Since $J^{-M}$, for $M\geqslant M_1$, and $S^{\infty}$ are homotopically equivalent by Lemma \ref{cinflema}, it follows from \cite[Corollary $2.11$]{AH} that the homology groups $H_{q}(J^{-M})$ and $H_{q}(S^{\infty})$ are isomorphic, for all $q\in\Z$. Since $S^{\infty}$ is also contractible in $X$ (see Benyamini {\it et al.} \cite{BS1}), it follows that  $\widetilde{H}_{*}(J^{-M})$ has the same reduced homology type of a point; namely,
\begin{equation}
    \widetilde{H}_{q}(J^{-M})\cong 0,\quad\mbox{ for all }q\in\Z. 
    \label{hgr1} 
\end{equation}

Next, we proceed to compute the critical groups of $J$ at infinity. Assume that $M\geqslant M_1$, and consider the topological pair $(X,J^{-M})$ and its long exact sequence of reduced homology groups given by 
\begin{equation}
\ldots{\rightarrow}{\widetilde H}_{q}(J^{-M})\stackrel{i_{*}}{\rightarrow}{\widetilde H}_{q}(X)\stackrel{j_{*}}\rightarrow {\widetilde H}_{q}(X,J^{-M})\stackrel{\partial_{*}}\rightarrow {\widetilde H}_{q-1}(J^{-M})\stackrel{i_{*}}\rightarrow\ldots,
\label{longseq}
\end{equation}
where $i_{*}$ and $j_{*}$ are the induced homomorphisms of the inclusion maps
$$i:J^{-M}\rightarrow X,\quad j:(X,\emptyset)\rightarrow (X,J^{-M}),$$ respectively.

Thus, since the sequence in (\ref{longseq}) is exact, it follows from (\ref{hgr1}) and the fact that $X$ is also contractible that 
\begin{equation}
C_{q}(J,\infty)=\widetilde{H}_q(X,J^{-M})\cong 0,\quad\mbox{ for }q \in \Z.
\label{cinftyg}
\end{equation}

\subsection{Existence of a Second Nontrivial Solution}
Assume that $J\in C^{2}(X,\R)$ and $z_o$ is a critical point of $J$. The {\it Morse index} of $z_o$ is defined as the supremum of the dimensions of the vector spaces of $X$ on which $J^{\prime\prime}(z_o)$ is negative definite. The {\it nullity} of $z_o$ is defined as the dimension of the $\ker J^{\prime\prime}(z_o).$ We will also denote $\mu^{*}=\mu_o+\nu_o$ as the {\it large Morse index}.

Next, we present a proposition about the computation of the critical groups at degenerate critical points. This result follows as a consequence of the splitting lemma and a result due to Gromoll and Meyer in \cite{Gr}.  We present the result as found in Cingolani and Vannella \cite{Cing2}.

\begin{theorem}\cite[Proposition $2.5$]{Cing2}
Suppose $X$ is a Hilbert space and $J\in C^{2}(X,\R)$. Let $z_o$ be an isolated critical point of $J$ with Morse index $\mu_o$ and large Morse index $\mu^\ast.$ Suppose $J^{\prime\prime}(z_o)$ is a Fredholm operator and let $V$ be the kernel of $J^{\prime\prime}(z_o)$. 
\begin{itemize}
    \item [(a)] If $z_o$ is a local minimum of $\widetilde{J}=J|_{V}$ then 
    $$C_{q}(J,z_o)\cong \delta_{q,\mu_o}\F,\quad\mbox{ for }q\in\Z.$$

    \item [(b)] If $z_o$ is a local maximum of $\widetilde{J}=J|_{V}$ then 
      $$C_{q}(J,z_o)\cong \delta_{q,\mu^{\ast}}\F,\quad\mbox{ for }q\in\Z.$$

    \item [(c)] If $z_o$ is neither a local maximum or a local minimum of $\widetilde{J}$, then 
    $$C_{q}(J,z_o)\cong 0,\quad\mbox{ for }q\not\in (\mu_0,\mu^\ast).$$
\end{itemize}
    \label{cingotheo}
\end{theorem}

Next, we present the main result of this section.

\begin{theorem}     Let $\Omega$ be a bounded and connected domain in $\R^{N}$, with $N\geqslant 3$, and suppose that all the hypotheses of Theorem  \ref{maintheo444} hold true.  In addition, assume 
condition $(H_4)$.
Let $z_o$ be the critical point of $J$ found in Theorem \ref{maintheo444}.   
 \begin{itemize}
     \item [(i)] If $z_o$ is a nondegenerate critical point of $J$ with Morse index $\mu_o$, then the system (\ref{mainsys}) has at least two nontrivial solutions provided that $q_o\ne \mu_o-1,$ where $q_o=\dim X^{-}.$
     \item [(ii)] If $z_o$ is a degenerate critical point of $J$ with Morse index $\mu_o$ and nullity $\nu_0$, then the system (\ref{mainsys}) has at least two nontrivial solutions provided that 
     \begin{equation} 
     q_o+1 \not\in[\mu_o,\mu_o+\nu_o].
     \label{qcond}
     \end{equation}
     where $q_o=\dim X^{-}.$ 
 \end{itemize}
\label{maintheo3}
\end{theorem}

\begin{proof}[Proof:]
Assume, by a way of contradiction, that the critical set $\mathcal{K}$ of $J$ has only two critical points; namely,  $\mathcal{K}=\{0,z_0\},$  where $z_o$ is the critical point found in Theorem \ref{maintheo4}. 
Set $c_0=J(z_0)$. Notice that $c_o>0$ as a result of Lemma \ref{lemmalink1}. 

Then, 
 by virtue of the definition of $M$ in Lemma \ref{cinflema}, we can choose $b_{1},b_{2},b_{3}$ such that
$$b_1<-M  <0 < b_2< c_{0} < b_{3},$$
where $M\geqslant M_1$. 
Thus, since $J$ satisfies the PS--condition, we obtain 
\begin{equation}
H_{q}(J^{b_{2}},J^{b_{1}})\cong C_{q}(J,0),\quad\mbox{ for }q\in\Z.
\label{cond1}
\end{equation}
Similarly, we also have 
\begin{equation}
H_{q}(J^{b_{3}},J^{b_{2}})\cong C_{q}(J,z_o),\quad\mbox{ for }q\in\Z.
\label{cond2}
\end{equation}
Consider the triple $(J^{b_{1}},J^{b_{2}},J^{b_{3}})$ and its long exact sequence given by
\begin{equation}
\ldots{\rightarrow}H_{k+1}(J^{b_{2}},J^{b_{1}})\stackrel{i_{*}}{\rightarrow}H_{k+1}(J^{b_{3}},J^{b_{1}})\stackrel{j_{*}}\rightarrow H_{k+1}(J^{b_{3}},J^{b_{2}})\stackrel{\partial_{*}}\rightarrow H_{k}(J^{b_{2}},J^{b_{1}})\stackrel{i_{*}}\rightarrow\ldots
\label{longseq2}
\end{equation}
(see \cite[Section $2.1$, page 118]{AH}).

By virtue of (\ref{cinftyg}), we also have 
\begin{equation}
H_{q}(J^{b_{3}},J^{b_{1}})=C_{q}(J,\infty)\cong \delta_{q,0}\F, \quad\mbox{ for } q\in\Z.
\label{lrr1}
\end{equation}

Thus, by virtue of (\ref{cond1}), (\ref{cond2}), and (\ref{lrr1}), we can write the sequence (\ref{longseq2}) as  
\begin{equation}
\ldots{\rightarrow}C_{k+1}(J,0)\stackrel{i_{*}}{\rightarrow}0\stackrel{j_{*}}\rightarrow C_{k+1}(J,z_o) \stackrel{\partial_{*}}\rightarrow C_{k}(J,0)\stackrel{i_{*}}\rightarrow 0\stackrel{j_{*}}\rightarrow\ldots,
\label{longseq22}
\end{equation}
for $k\geqslant 1$.

Using the fact that the sequence (\ref{longseq22}) is exact, we claim that $\partial_{*}$ is an isomorphism. In fact, $\ker(\partial_{*})=Im(j_{*})=\{0\}$ implies that $\partial_{*}$ is one-to-one and $Im(\partial_*)=\ker (i_{*})$ implies that $\partial_*$ is onto because $Im(i_*)=\{0\}$. 

It follows from (\ref{cgrouporigin}), (\ref{cond2}), and the fact that $\partial_*$ is an isomorphism in (\ref{longseq22}), that 
\begin{equation}
   \delta_{k,q_o}\F \cong C_{k}(J,0)\cong C_{k+1}(J,z_o),\quad\mbox{ for }k\in\Z,
   \label{ll1}
\end{equation}
where $q_o=\dim X^{-}.$ 

Next, we analyze the the cases when $z_o$ is a nondegenerate critical point or a degenerate critical point of the functional $J$ defined in (\ref{func0}). We will show that we obtain a contradiction in each of these cases. 

For the case when $z_o$ is a nondegenerate critical point for $J$ with Morse index $\mu_o$, we invoke the result in \cite[Corollary $8.3$]{MW1} to get 
$$C_{k}(J,z_o)\cong \delta_{k,\mu_o}\F,\quad\mbox{ for }k\in\Z,$$
where $\mu_o$ is the Morse index of $z_o.$  

Then, by virtue of (\ref{ll1}) and the previous isomorphism, we get 
\begin{equation}
    \delta_{k,q_o}\F \cong C_{k+1}(J,z_o)\cong \delta_{k+1,\mu_o}\F\cong \delta_{k,\mu_o-1}\F,\quad\mbox{ for }k\in\Z.
    \label{shift66}
\end{equation}

Set $k=q_o$ in (\ref{shift66}) to conclude that $q_o=\mu_o-1$ which is a contradiction since $q_o\ne \mu_o-1$ by assumption.

Next, let $z_o$ be a degenerate critical point for $J$ with Morse index $\mu_0$ and nullity $\nu_0$. Then, it follows from the shifting theorem \cite[Theorem $8.4$] {MW1}  that 
\begin{equation}
    C_{k}(J,z_o)\cong C_{k-\mu_o}(\widetilde{J},z_o),\quad\mbox{ for }k\in\Z,
    \label{shift1}
\end{equation}
where $\widetilde{J}:B_{\rho}(z_o)\rightarrow \R$ is a $C^2$ function, $B_{\rho}(z_o)$ is a small ball in $V_o:=\ker J^{\prime\prime}(z_o)$ centered at $z_o$ and radius $\rho>0$, $\nu_o=\dim V_o.$

In what follows, we analyze the three possible cases for the critical point $z_o$ as in Theorem \ref{cingotheo} and show that we reach a contradiction in each one of them. 
 
If $z_o$ was a local minimum for $\widetilde{J}$, then, by virtue of (\ref{ll1}), (\ref{shift1}), and \cite[Example $1$, page $33$]{KC}, we would have 
\begin{equation}
    \delta_{k,q_o}\F\cong C_{k+1}(J,z_o)\cong C_{k+1-\mu_o}(\widetilde{J},z_o)\cong \delta_{k+1-\mu_o,0}\F\cong \delta_{k+1,\mu_o}\F, 
 \nonumber 
\end{equation}
for all $k\in\Z$. Set $k=q_o$ in the previous equation to conclude that $\F\cong 0$ since $q_o+1\ne \mu_o$ by assumption. This is a contradiction.

If $z_o$ was a local maximum for $\widetilde{J}$, then by virtue of (\ref{ll1}), (\ref{shift1}), and  \cite[Example $2$, page $33$]{KC}, 
\begin{equation}
   \delta_{k,q_o}\F\cong  C_{k+1}(J,z_o)\cong  C_{k+1-\mu_o}(\widetilde{J},z_o)\cong \delta_{k+1-\mu_o,\nu_o}\F \cong \delta_{k+1,\mu_o+\nu_o}\F,
   \label{lasteq}
\end{equation}
for all $k\in\Z.$ Observe that, if we set $k=q_o$ in (\ref{lasteq}), then we conclude that $\F\cong 0$ since $q_o+1\ne\mu_o+\nu_o$ by virtue of (\ref{qcond}).

 Finally, if $z_o$ is neither a local maximum nor a local minimum for $\widetilde{J}$, it follows from Theorem \ref{cingotheo}-(c) that  
 \begin{equation} 
 C_{q}(J,z_o)\cong 0,\quad\mbox{ for } q\not\in(\mu_o,\mu_o+\nu_o), 
 \label{ss4}
 \end{equation}
 
 Hence, by virtue of (\ref{ll1}) and  (\ref{ss4}), we get
 \begin{equation}
     \delta_{k,q_o}\F\cong C_{k+1}(J,z_o)\cong 0,\quad\mbox{ for }k+1\not\in (\mu_o,\mu_o+\nu_o).
     \label{cc11}
 \end{equation}

Set $k=q_o$ in (\ref{cc11}) to get $\F\cong 0,$ which is a contradiction. 

Therefore, the critical set $\mathcal{K}$ must have at least three critical points. This implies that problem (\ref{mainsys}) has at least two nontrivial solutions, and this concludes the proof of the theorem. 

\end{proof}

\bibliographystyle{plain}
\bibliography{main}

\end{document}